\documentclass[11pt, reqno]{amsart}
\usepackage{amssymb}
\usepackage{verbatim}
\usepackage[vcentermath]{youngtab}
\usepackage{float}
\usepackage{hyperref}
\usepackage{tikz}
\usepackage{hyperref}
\usepackage{multicol}
\usepackage{wrapfig}
\usepackage{multicol}
\usepackage[left=1.1in, right=1.1in, top=1in, bottom=.9in]{geometry}

\usepackage{times}
\usepackage[T1]{fontenc}
\usepackage{mathrsfs}
\usepackage{epsfig}
\usepackage{color}
\usepackage{array} 
\newcolumntype{L}{>{$}l<{$}}
\newcolumntype{R}{>{$}r<{$}}

\usepackage{enumitem}

\newtheorem{thm}{Theorem}[section]
\newtheorem{lemma}[thm]{Lemma}

\newtheorem{cor}[thm]{Corollary}
\newtheorem{prop}[thm]{Proposition}

\theoremstyle{remark}

\theoremstyle{definition}

\newtheorem{remark}[thm]{Remark}
\newtheorem{caution}[thm]{Caution}

\newtheorem{examplex}[thm]{Example}

\newenvironment{example}
  {\pushQED{\qed}\examplex}
  {\popQED\endexamplex}

\def\ZZ{\mathbb{Z}}
\def\QQ{\mathbb{Q}}

\def\NN{\mathbb{N}}

\def\AA{\mathbb{A}}

\def\FF{\mathbb{F}}

\def\RR{\mathbb{R}}
\def\CC{\mathbb{C}}

\def\C{\mathcal{C}}
\def\U{\mathcal{U}}

\def\H{\mathcal{H}}

\def\PConf{\mathrm{PConf}}
\def\wt{\widetilde}

\def\h{\widehat}
\def\ve{\varepsilon}

\def\U{\mathcal{U}}

\def\s{\hspace{.1em}}

\def\multi#1#2{\ensuremath{\left(\kern-.3em\left(\genfrac{}{}{0pt}{}{#1}{#2}\right)\kern-.3em\right)}}

\usepackage{empheq}

\numberwithin{equation}{section}

\begin{document}

\title[Cyclotomic Factors of Necklace Polynomials]{Cyclotomic Factors of Necklace Polynomials}

\author{Trevor Hyde}
\address{Dept. of Mathematics\\
University of Chicago \\
Chicago, IL 60637\\
}
\email{tghyde@uchicago.edu}

\maketitle

\begin{abstract}
    We observe that the necklace polynomials $M_d(x) = \frac{1}{d}\sum_{e\mid d}\mu(e)x^{d/e}$ are highly reducible over $\QQ$ with many cyclotomic factors. Furthermore, the sequence $\Phi_d(x) - 1$ of shifted cyclotomic polynomials exhibits a qualitatively similar phenomenon, and it is often the case that $M_d(x)$ and $\Phi_d(x) - 1$ have many common cyclotomic factors. We explain these cyclotomic factors of $M_d(x)$ and $\Phi_d(x) - 1$ in terms of what we call the \emph{$d$th necklace operator}. Finally, we show how these cyclotomic factors correspond to certain hyperplane arrangements in finite abelian groups.
\end{abstract}

\section{Introduction}
\label{sec intro}
The \emph{$d$th necklace polynomial} $M_d(x)$, for positive integral $d$, is defined by
\[
    M_d(x) := \frac{1}{d}\sum_{e\mid d}\mu(e) x^{d/e},
\]
where $\mu$ is the number theoretic M\"obius function and the sum is over all divisors $e$ of $d$. Necklace polynomials arise naturally in number theory, combinatorics, dynamics, geometry, representation theory, and algebra. For example, if $q$ is a prime power and $\FF_q$ is a finite field with $q$ elements, then $M_d(q)$ is the number of $\FF_q$-irreducible monic polynomials of degree $d$ in $\FF_q[x]$; if $k\geq 1$ is a natural number, then $M_d(k)$ is the number of aperiodic necklaces comprised of $d$ beads chosen from among $k$ colors.

We begin with the empirical observation that necklace polynomials are highly reducible over $\QQ$. For example, if $d = 105$, then
\begin{align}
\label{eqn 105 factor}
    M_{105}(x) &= \tfrac{1}{105}(x^{105} - x^{35} - x^{21} - x^{15} + x^7 + x^5 + x^3 - x)\nonumber\\
    &= e(x)(x^4 + 1)(x^2 - x + 1)(x^2 + 1)(x^2 + x + 1)(x + 1)(x - 1)x,
\end{align}
where $e(x)\in \QQ[x]$ is an irreducible polynomial of degree 92. With only two exceptions, the irreducible factors of $M_{105}(x)$ are cyclotomic polynomials. Recall that the \emph{$m$th cyclotomic polynomial} $\Phi_m(x)$ is the $\QQ$-minimal polynomial of a primitive $m$th root of unity. With this notation \eqref{eqn 105 factor} may be expressed as
\[
    M_{105}(x) = e(x)\cdot \Phi_8\cdot \Phi_6\cdot \Phi_4\cdot \Phi_3\cdot \Phi_2\cdot \Phi_1\cdot x.
\]
Here are several more examples: There are irreducible, non-cyclotomic polynomials $f(x), g(x), h(x) \in \QQ[x]$ with degrees $148, 212,$ and  $708$, respectively, such that
\begin{align*}
    M_{165}(x) &= \tfrac{1}{165}(x^{165} - x^{55} - x^{33} - x^{15} + x^{11} + x^{5} + x^{3} - x)\\
    &= f(x)\cdot \Phi_{12} \cdot \Phi_{10} \cdot \Phi_5 \cdot \Phi_4 \cdot \Phi_2 \cdot \Phi_1 \cdot x\\
    M_{231}(x) &= \tfrac{1}{231}(x^{231} - x^{77} - x^{33} - x^{21} + x^{11} + x^{7} + x^{3} - x)\\
    &= g(x) \cdot \Phi_{10} \cdot \Phi_8 \cdot \Phi_6 \cdot \Phi_5 \cdot \Phi_3 \cdot \Phi_2 \cdot \Phi_1 \cdot x\\
    M_{741}(x) &= \tfrac{1}{741}(x^{741} - x^{247} - x^{57} - x^{39} + x^{19} + x^{13} + x^3 - x)\\
    &= h(x)\cdot \Phi_{20}\cdot \Phi_{18}\cdot \Phi_{12}\cdot \Phi_9\cdot \Phi_6\cdot \Phi_4\cdot \Phi_3\cdot \Phi_2\cdot \Phi_1\cdot x.\qedhere
\end{align*}
Since $M_d(x)$ has rational coefficients, $\Phi_m(x)$ dividing $M_d(x)$ is equivalent to $M_d(\zeta_m) = 0$ for some primitive $m$th root of unity $\zeta_m$. The plot below shows all pairs $(d,m)$ with $1\leq d, m \leq 1000$ such that $M_d(\zeta_m) = 0$.
\begin{center}
    \includegraphics[scale=.7]{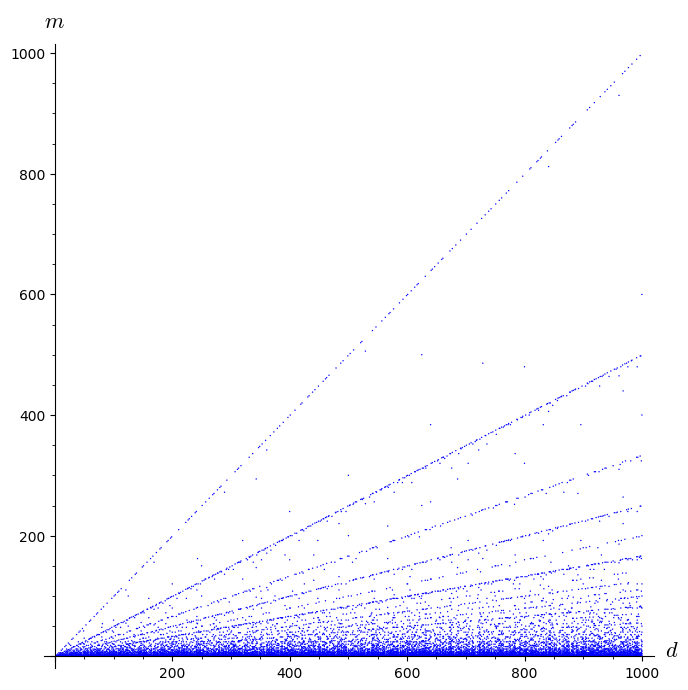}
\end{center}
This plot suggests that the preponderance of cyclotomic factors of $M_d(x)$ observed above is not isolated to special values of $d$, but rather that it occurs to some extent for all $d$. The primary objectives of this paper are to explain why necklace polynomials have so many cyclotomic factors and to characterize the pairs of integers $(d,m)$ for which $M_d(\zeta_m) = 0$.

A strikingly similar phenomenon occurs for the seemingly unrelated sequence $\Phi_d(x) - 1$ of shifted cyclotomic polynomials. For example,
\begin{align*}
    \Phi_{105}(x) - 1 &= \wt{e}(x) \cdot \Phi_8 \cdot \Phi_6 \cdot \Phi_4 \cdot \Phi_3 \cdot \Phi_2 \cdot \Phi_1\cdot x\\
    \Phi_{165}(x) - 1 &= \wt{f}(x) \cdot \Phi_{10} \cdot \Phi_5 \cdot \Phi_4 \cdot \Phi_2 \cdot \Phi_1\cdot x\\
    \Phi_{231}(x) - 1 &= \wt{g}(x) \cdot \Phi_{12} \cdot \Phi_{10} \cdot \Phi_6 \cdot \Phi_5 \cdot \Phi_4 \cdot \Phi_3 \cdot \Phi_2 \cdot \Phi_1 \cdot x\\
    \Phi_{741}(x) - 1 &= \wt{h}(x) \cdot \Phi_{18} \cdot \Phi_{12} \cdot \Phi_{9} \cdot \Phi_{6} \cdot \Phi_{4} \cdot \Phi_{3} \cdot \Phi_{2} \cdot \Phi_{1} \cdot x,
\end{align*}
where $\wt{e}(x), \wt{f}(x), \wt{g}(x), \wt{h}(x) \in \ZZ[x]$ are irreducible, non-cyclotomic polynomials with degrees 35, 67, 99, and 407, respectively. Note that $\Phi_m(x)$ dividing $\Phi_d(x) - 1$ is equivalent to $\Phi_d(\zeta_m) = 1$ for a primitive $m$th root of unity $\zeta_m$.

Comparing the factorizations of $M_d(x)$ and $\Phi_d(x) - 1$ in the examples above we see there is a considerable overlap in their cyclotomic factors. The table below illustrates that this is a common occurrence. For each $2 \leq d \leq 43$, we list all $m$ for which $\Phi_m(x)$ divides both $M_d(x)$ and $\Phi_d(x) - 1$ in plain text, and all $m$ for which $\Phi_m(x)$ divides $M_d(x)$ but not $\Phi_d(x) - 1$ in bold. For $d$ in this range, there are no $m$ for which $\Phi_m(x)$ divides $\Phi_d(x) - 1$ but not $M_d(x)$; the first time this occurs is with $d = 231$ and $m = 4$.

The secondary objectives of this paper are to explain why this qualitatively similar cyclotomic factor phenomenon occurs for the shifted cyclotomic polynomials $\Phi_d(x) - 1$, explain how these factors are related to the factors of $M_d(x)$, and to characterize those pairs of integers $(d,m)$ for which $\Phi_d(\zeta_m) = 1$.

\begin{center}
\resizebox{.7\textwidth}{!}{
    \begin{tabular}{|c|l||c|l||c|l|}
    \hline
        $d$ & $m$ & $d$ & $m$ & $d$ & $m$\\
    \hline
         $2$ & $\mathbf{1} $                            & $16 $ & $\mathbf{1}, \mathbf{2}, \mathbf{4}, \mathbf{8} $ & $30 $ & $1,2,4,\mathbf{6} $\\
         $3$ & $\mathbf{1}, 2 $                         & $17 $ & $\mathbf{1}, 2, 4, 8, 16 $                        & $31 $ & $\mathbf{1},2,3,5,6,10,15,30 $\\
         $4$ & $\mathbf{1}, \mathbf{2} $                & $18 $ & $1, \mathbf{2}, 3, \mathbf{6} $                   & $32 $ & $\mathbf{1},\mathbf{2},\mathbf{4},\mathbf{8},\mathbf{16} $\\
         $5$ & $\mathbf{1}, 2, 4 $                      & $19 $ & $\mathbf{1}, 2, 3, 6, 9, 18 $                     & $33 $ & $1,2,5,10 $\\
         $6$ & $1, \mathbf{2} $                         & $20 $ & $1, 2, \mathbf{4}, 8, \mathbf{12} $               & $34 $ & $1,2,4,\mathbf{6},8,16 $\\
         $7$ & $\mathbf{1}, 2, 3, 6 $                   & $21 $ & $1, 2, \mathbf{3}, 6, \mathbf{8} $                & $35 $ & $1,2,3,4,6 $\\
         $8$ & $\mathbf{1}, \mathbf{2}, \mathbf{4} $    & $22 $ & $1, \mathbf{2}, 5, \mathbf{6}, 10 $               & $36 $ & $1,2,3,\mathbf{4},6,\mathbf{12} $\\
         $9$ & $\mathbf{1}, 2, \mathbf{3}, 6 $          & $23 $ & $\mathbf{1}, 2, 11, 22 $                          & $37 $ & $\mathbf{1},2,3,4,6,9,12,18,36 $\\
         $10$ & $1, \mathbf{2}, 4, \mathbf{6} $         & $24 $ & $1, 2, 4, \mathbf{8} $                            & $38 $ & $1,\mathbf{2},3,6,9,18 $\\
         $11$ & $\mathbf{1}, 2, 5, 10 $                 & $25 $ & $\mathbf{1}, 2, 4, \mathbf{5}, 10, 20 $           & $39 $ & $1,2,\mathbf{3},4,6,12 $\\
         $12$ & $1, 2, \mathbf{4} $                     & $26 $ & $1, \mathbf{2}, 3, 4, 6, 12 $                     & $40 $ & $1,2,4,\mathbf{8},16,\mathbf{24} $\\
         $13$ & $\mathbf{1}, 2, 3, 4, 6, 12 $           & $27 $ & $\mathbf{1}, 2, \mathbf{3}, 6, \mathbf{9}, 18 $   & $41 $ & $\mathbf{1},2,4,5,8,10,20,40 $\\
         $14$ & $1, \mathbf{2}, 3, 6 $                  & $28 $ & $1, 2, 3, \mathbf{4}, 6, 12 $                     & $42 $ & $1,2,3,\mathbf{6} $\\
         $15$ & $1, 2, 4 $                              & $29 $ & $\mathbf{1}, 2, 4, 7, 14, 28 $                    & $43 $ & $\mathbf{1},2,3,6,7,14,21,42 $\\
    \hline
    \end{tabular}
    }
\end{center}
\vspace{.2in}

\begin{wrapfigure}{r}{.5\textwidth}
\begin{center}
    \includegraphics[scale=.165]{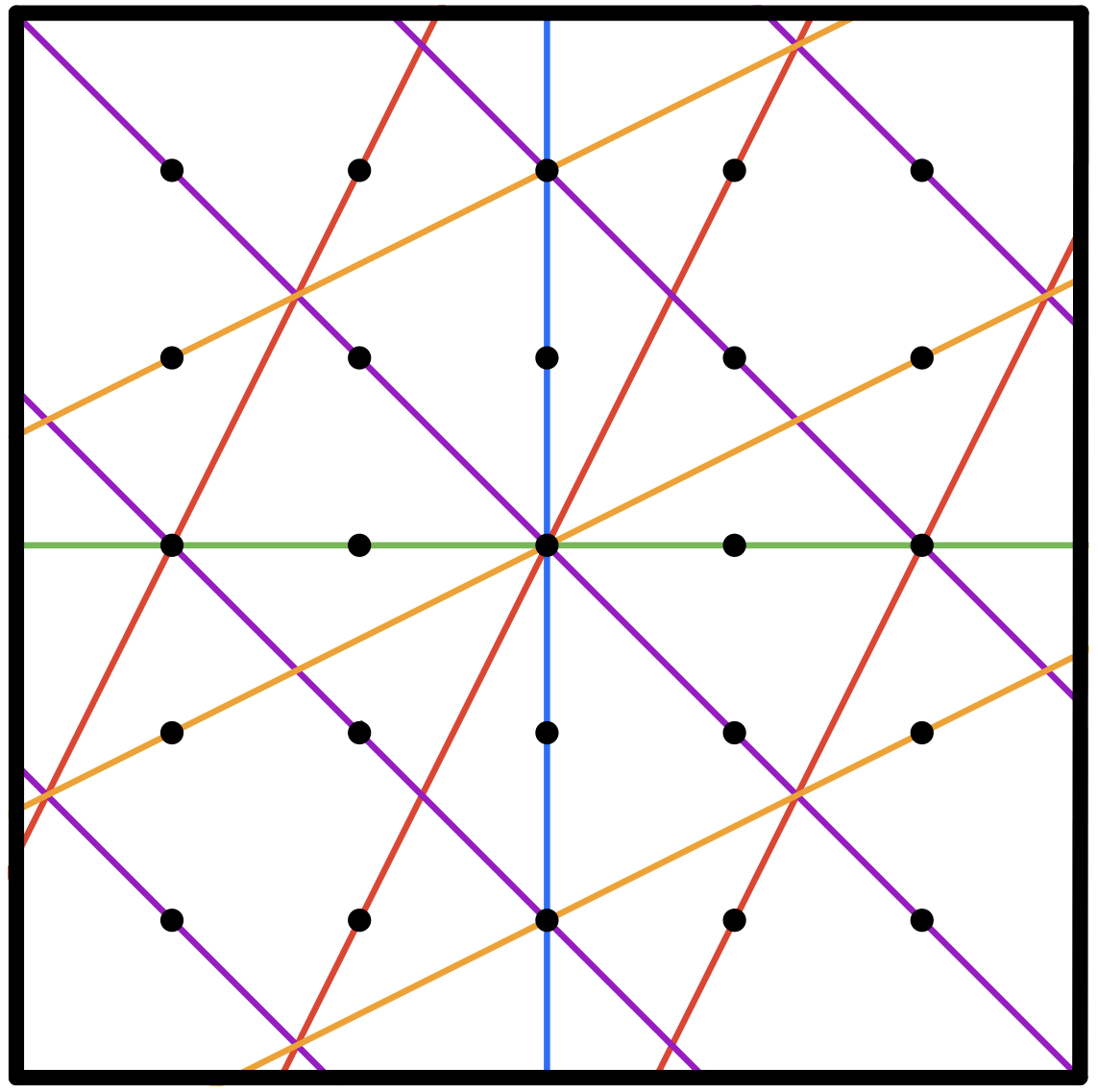}
\end{center}
\end{wrapfigure}

We explain the cyclotomic factors of necklace polynomials $M_d(x)$ and shifted cyclotomic polynomials $\Phi_d(x) - 1$ using the representation theory of finite abelian groups. We trace this phenomenon in both cases to a common source, which we call the \emph{necklace operators}, and show how these operators account for the common cyclotomic factors of $M_d(x)$ and $\Phi_d(x) - 1$. Our analysis reveals a surprising connection between these unexpected cyclotomic factors and arrangements of hyperplanes in finite abelian groups. For example, we will explain how the arrangement of lines covering $\ZZ/(4) \times \ZZ/(4)$ pictured to the right corresponds to the fact that $M_d(\zeta_{65}) = 0$ and $\Phi_d(\zeta_{65}) = 1$ with $d = 9372603371$ (see Example \ref{ex lines}.) Our terminology and explicit results are detailed in the following section.

\subsection{Results}

Our first result relates the identities $M_d(\zeta_m) = 0$ and $\Phi_d(\zeta_m) = 1$ and hyperplane arrangements in the group of Dirichlet characters of modulus $m$. Let $\s\U_m := (\ZZ/(m))^\times$ denote the multiplicative group of integers modulo $m$ and let $\h{\U}_m := \mathrm{Hom}(\s\U_m,\CC^\times)$ be the group of \emph{Dirichlet characters of modulus $m$}. Each unit $q \in \U_m$ determines a homomorphism from $\h{\U}_m$ to $\CC^\times$ by $\chi \mapsto \chi(q)$; let $\H_q \subseteq \h{\U}_m$ denote the kernel of this map. We call $\H_q$ the \emph{hyperplane associated to $q$},
\[
    \H_q := \{\chi \in \h{\U}_m : \chi(q) = 1\}.
\]
Note that with a choice of coordinates for the group $\,\h{\U}_m$---by which we mean some isomorphism between $\h{\U}_m$ and a product of cyclic groups $\ZZ/(n)$---$\H_q$ may be expressed as the vanishing set of an integral linear form, hence the hyperplane terminology (see Remark \ref{ex hyperplanes}.)

\begin{thm}
\label{thm simple intro}
Let $d, m > 1$ be coprime integers. If $\s\h{\U}_m \subseteq \bigcup_{p \mid d} \H_p$, then $x^m - 1$ divides $M_d(x)$ and $\frac{x^m - 1}{x-1}$ divides $\Phi_d(x) - 1$.
\end{thm}

In other words, if the group $\h{\U}_m$ of Dirichlet characters of modulus $m$ is covered by the arrangement of hyperplanes $\{\H_p : p \mid d \text{ is prime}\}$, then $M_d(\zeta_m^k) = 0$ for all $k\geq 0$ and $\Phi_d(\zeta_m^k) = 1$ for all $k \not\equiv 0 \bmod m$.

\begin{remark}
Theorem \ref{thm simple intro} avoids addressing $\Phi_d(1)$, but it is well-known that $\Phi_d(1) = 1$ whenever $d$ is divisible by at least two distinct primes and that $\Phi_{p^r}(1) = p$ for any prime $p$ and $r\geq 1$.
\end{remark}

Theorem \ref{thm simple intro} shows that hyperplane arrangements covering $\h{\U}_m$ provide one source of common cyclotomic factors of $M_d(x)$ and $\Phi_d(x) - 1$, and that these factors have the property that if $\Phi_m(x)$ is a factor, so is $\Phi_n(x)$ for all $n>1$ dividing $m$. Theorem \ref{thm simple intro} empirically accounts for the majority of such common cyclotomic factors. For example, with $1 \leq d \leq 1000$, Theorem \ref{thm simple intro} accounts for all common cyclotomic factors of $M_d(x)$ and $\Phi_d(x) - 1$; for about $88.9\%$ of the cyclotomic factors of $M_d(x)$; and for about $99.7\%$ of the cyclotomic factors of $\Phi_d(x) - 1$.

\begin{example}
\label{ex m = 24 intro}
We illustrate Theorem \ref{thm simple intro} in the case $m = 24$. The Dirichlet characters $\h{\U}_{24}$ form a 3 dimensional $\FF_2$-vector space. Note that $\U_{24}$ is generated by 13, 17, and 19. Identifying $\U_{24}$ with the dual of $\h{\U}_{24}$ we can choose coordinates $\rho: \U_{24} \rightarrow \h{\FF}_2^3$ such that $\rho(13) = x$, $\rho(17) = y,$ and  $\rho(19) = z$. The pencil of planes containing the line $\langle (1,1,1) \rangle$ covers all of $\FF_2^3 \cong \h{\U}_{24}$ and consists of
\[
    \H_{13\cdot 17} : x + y = 0,\hspace{.4in}
    \H_{13\cdot 19} : x + z = 0,\hspace{.4in}
    \H_{17\cdot 19} : y + z = 0.
\]
Since
\[
    13\cdot 17 \equiv 5 \bmod 24, \hspace{.35in} 13\cdot 19 \equiv 7 \bmod 24, \hspace{.35in} 17\cdot 19 \equiv 11 \bmod 24,
\]
it follows from Theorem \ref{thm simple intro} with $d = 385 = 5 \cdot 7 \cdot 11$ that $x^{24} - 1$ divides $M_{385}(x)$ and $\frac{x^{24}-1}{x - 1}$ divides $\Phi_{385}(x) - 1$.
\begin{center}
    \includegraphics[scale=.2]{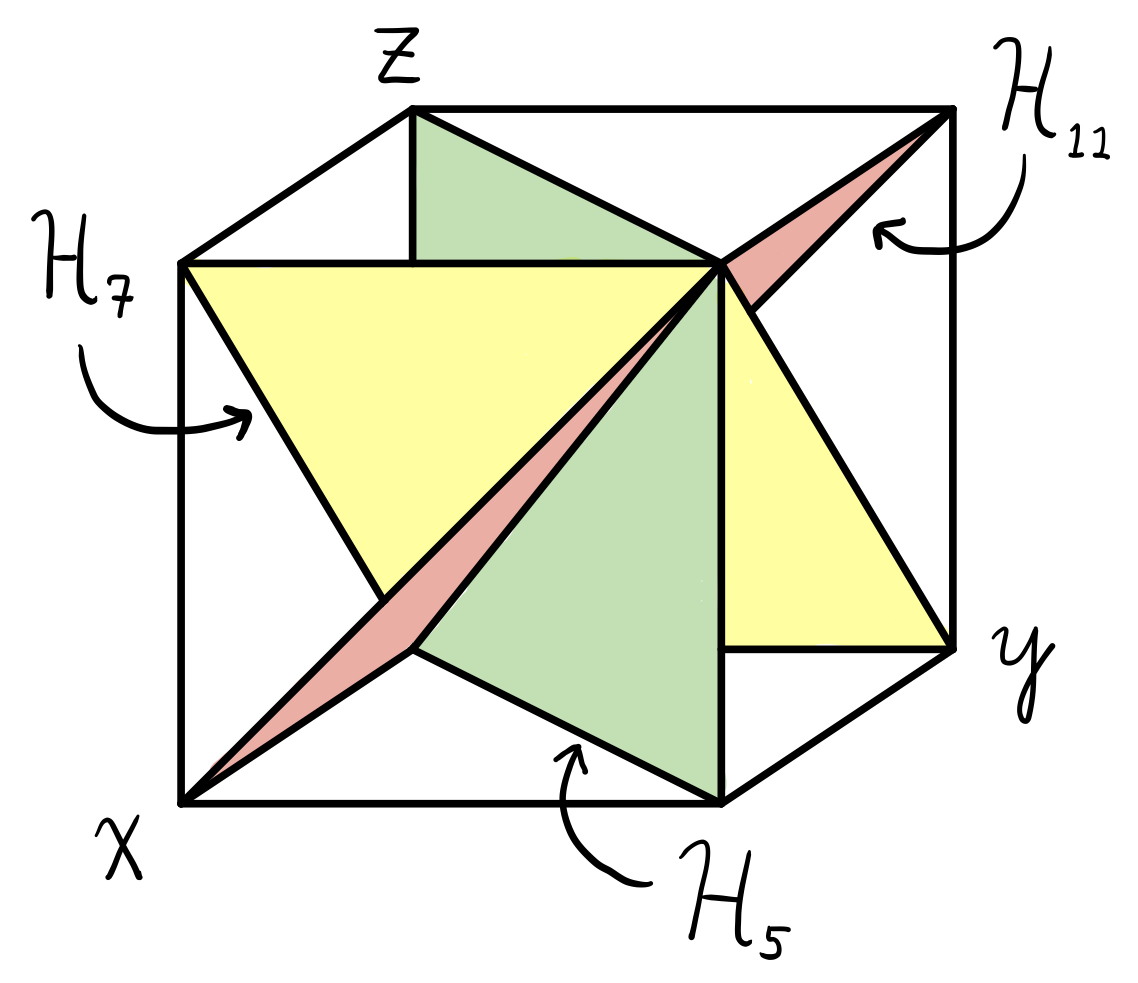}
\end{center}
\end{example}

\begin{example}
\label{ex triv hyp}
Let $d, m \geq 1$ and suppose that $d$ is divisible by some prime $p$ such that $p \equiv 1 \bmod m$. In this case, $\H_p = \H_1 = \h{\U}_m$ is the degenerate hyperplane, namely the entire group (recall that $\H_d$ is the kernel of the evaluation map $\chi \mapsto \chi(d)$ for $\chi$ a Dirichlet character of modulus $m$.) Hence the arrangement $\{\H_p \subseteq \h{\U}_m : p \mid d \text{ is prime}\}$ trivially covers $\h{\U}_m$. Thus Theorem \ref{thm simple intro} implies that $M_d(\zeta_m) = 0$ and $\Phi_d(\zeta_m) = 1$ whenever $d$ is divisible by a prime $p$ such that $p \equiv 1 \bmod m$. In particular, with $d$ fixed, this holds for $m = p - 1$ if $\gcd(d, p - 1) = 1$. This explains why cyclotomic factors of $M_d(x)$ and $\Phi_d(x) - 1$ are so prevalent: each such prime $p$ dividing $d$ contributes a factor of $\frac{x^{p-1}-1}{x - 1}$ to both polynomials.
\end{example}

Our second result highlights the structure of the pairs $(d,m)$ with $m$ fixed for which $M_d(\zeta_m) = 0$ or $\Phi_d(\zeta_m) = 1$.
\begin{thm}
\label{thm primewise intro}
Let $d, e, m \geq 1$.
\begin{enumerate}
    \item If $M_d(\zeta_m) = 0$ and $e$ is coprime to $m$, then $M_{de}(\zeta_m) = 0$.
    
    \item If $d$ and $e$ are coprime to $m$ and if we have an equality of sets of residue classes
    \[
        \{p \bmod m : p \mid d \text{ is prime}\} = \{q \bmod m : q \mid e \text{ is prime}\},
    \]
    then $M_d(\zeta_m) = 0$ if and only if $M_e(\zeta_m) = 0$.
\end{enumerate}
Likewise both assertions hold with $M_d(\zeta_m) = 0$ replaced by $\Phi_d(\zeta_m) = 1$.
\end{thm}

\begin{example}
In our examples above we saw that $M_{231}(\zeta_8) = 0$. Thus Theorem \ref{thm primewise intro}.1 implies that $M_{231e}(\zeta_8) = 0$ for all odd $e$ and Theorem \ref{thm primewise intro}.2 implies that $M_d(\zeta_{8}) = 0$ for $d$ any product of odd primes with at least one congruent to each of $3, 7 \bmod 8$, including, for instance, $M_{21}(\zeta_{8}) = M_{77}(\zeta_{8}) = 0$.
\end{example}

\begin{example}
A quick computation shows that $M_{10}(\zeta_6) = 0$ but $M_{20}(\zeta_6) \neq 0$. This example shows that the assumption that $e$ is coprime to $m$ is necessary in Theorem \ref{thm primewise intro}.1.
\end{example}

Theorem \ref{thm intro precise vanish} characterizes the pairs $(d,m)$ for which $M_d(\zeta_m) = 0$ or $\Phi_d(\zeta_m) = 1$, without the coprime restriction on $d$ and $m$, in terms of hyperplane arrangements covering certain prescribed subsets of $\h{\U}_m$. First, some set-up. If $n$ divides $m$, then there is a natural injective map $\h{\U}_n \rightarrow \h{\U}_m$ induced by the quotient $\U_m \rightarrow \U_n$. We use these maps to identify $\h{\U}_n$ with its image in $\h{\U}_m$ and say $\h{\U}_n \subseteq \h{\U}_m$. If $\chi \in \h{\U}_m$, then let $c_\chi$ be the smallest positive integer $n$ such that $\chi \in \h{\U}_{n}$. Finally, let $v_p$ denote the normalized $p$-adic valuation.

\begin{thm}
\label{thm intro precise vanish}
Let $d, e, f, m \geq 1$ be integers and let $m'$ be the product of all primes $p$ such that $v_p(m) = 1$. Suppose that
\begin{multicols}{2}
\begin{enumerate}[label = (\roman*)]
    \item $def$ is squarefree,
    \item $d$ is coprime to $m$,
    \item $e$ divides $m'$,
    \item $f$ divides $m/m'$.
\end{enumerate}
\end{multicols}
\begin{enumerate}
    \item If $\Sigma_{f,m} \subseteq \h{\U}_m$ is the set of all characters $\chi$ such that
\begin{enumerate}
    \item $v_p(c_\chi) = v_p(m)$ if $v_p(m) \geq 2$ and $v_p(f) = 0$, and
    \item $v_p(c_\chi) \geq v_p(m) - 1$ if $v_p(m) > 2$ and $v_p(f) = 1$,
\end{enumerate}
then $M_{def}(\zeta_m) = 0$ if and only if
\[
     \Sigma_{f,m}\subseteq \begin{cases} \bigcup_{p\mid d} \H_p & \text{if $2\nmid e$,} \\ \bigcup_{p\mid d} \H_p \cup \H_2^a & \text{if $2\mid e$,}\end{cases}
\]
where $\H_2^a \subseteq \h{\U}_m$ is the affine hyperplane $\H_2^a := \{\chi \in \h{\U}_m : \chi(2) = -1\}$.\\
    
\item If $m$ does not divide $def$, then $\Phi_{def}(\zeta_m) = 1$ if and only if
\begin{enumerate}
    \item $\displaystyle{\H_{-1} \subseteq \begin{cases} \bigcup_{p\mid md/e} \H_p & \text{if }3 \nmid e\\ \bigcup_{p\mid md/e} \H_p\cup \H_3^a & \text{if }3 \mid e, \end{cases}}$
    
    where $\H_3^a \subseteq \h{\U}_m$ is the affine hyperplane $\H_3^a := \{\chi \in \h{\U}_m : \chi(3) = -1\}$,\\
    \item $m$ divides $\varphi(def)$, and\\
    \item $\displaystyle{\sum_{a\mid def}\lfloor a/m \rfloor \equiv \frac{\varphi(def)}{m} \bmod 2}$.
\end{enumerate}

\end{enumerate}
\end{thm}

\begin{remark}
Several comments on Theorem \ref{thm intro precise vanish}.
\begin{enumerate}[leftmargin=*]
    \item Most of the subtlety in characterizing the pairs $(d,m)$ for which $M_d(\zeta_m) = 0$ or $\Phi_d(\zeta_m) = 1$ arises from common factors of $d$ and $m$. The essential point is that the identities $M_{def}(\zeta_m) = 0$ and $\Phi_{def}(\zeta_m) = 1$ correspond to a certain subset $\Sigma_{def} \subseteq \h{\U}_m$ of characters being covered by an arrangement of (affine) hyperplanes in $\h{\U}_m$.

    \item If $d \geq 1$, let $d_0$ be the product of all distinct primes dividing $d$ and let $e = d/d_0$. Then $dM_d(\zeta_m) = d_0M_{d_0}(\zeta_m^e)$ and $\Phi_d(\zeta_m) = \Phi_{d_0}(\zeta_m^e)$. Hence we lose no generality in Theorem \ref{thm intro precise vanish} by assuming that $def$ is squarefree.
    
    \item If $\chi \in \h{\U}_m$ is a character, then a common convention is to set $\chi(d) = 0$ whenever $d$ is not coprime to $m$. Our identification of $\h{\U}_n$ with its image in $\h{\U}_m$ induced by the quotient map $\U_m \rightarrow \U_n$ whenever $n$ divides $m$ suggests a slight natural variant on this convention which we find convenient: If $\chi \in \h{\U}_m$ has conductor $n$ and $d \in \ZZ$, then we set $\chi(d) = 0$ if $d$ is not coprime to $n$ and otherwise set $\chi(d)$ to its nonzero value on the residue class of $d$ modulo $n$. In particular, the characters on the affine hyperplane $\H_2^a$ defined in Theorem \ref{thm intro precise vanish} must all have conductor dividing $m/2$. See Caution \ref{caution}.
\end{enumerate}
\end{remark}

\begin{example}
\label{ex m = 8}
Theorem \ref{thm intro precise vanish} allows us to account for the cyclotomic factors of $M_d(x)$ not explained by Theorem \ref{thm simple intro}. For example, let $d = 21$. Then $M_{21}(x)$ factors as
\[
    M_{21}(x) = f(x)\cdot \Phi_8 \cdot \Phi_6 \cdot \Phi_3 \cdot \Phi_2 \cdot \Phi_1 \cdot x,
\]
where $f(x) = \tfrac{1}{10}(x^{10} - x^6 + x^4 + x^2 - 1)$ is irreducible and not cyclotomic. The factor $\Phi_8(x)$ cannot follow from Theorem \ref{thm simple intro} since $\Phi_4(x)$ does not divide $M_{21}(x)$, thus we turn to Theorem \ref{thm intro precise vanish}.

Using the notation of Theorem \ref{thm intro precise vanish}, we have $m = 8$ and $e = f = 1$. Since $8$ is a prime power with exponent at least 2, the set $\Sigma_{1,8}$ consists of the characters with conductor 8. There are two such characters $\chi$ determined by $\chi(3) = \pm 1$ and $\chi(5) = -1$. If $\chi(3) = 1$, then $\chi \in \H_3$ and if $\chi(3) = -1$, then 
\[
    \chi(7) = \chi(3)\chi(5) = (-1)^2 = 1,
\]
hence $\chi \in \H_7$. Thus $\Sigma_{1,8} \subseteq \H_3 \cup \H_7$ and Theorem \ref{thm intro precise vanish} implies that $M_{21}(\zeta_8) = 0$. We can visualize this situation with the following diagram: we choose coordinates for $\h{\U}_8 \cong \ZZ/(2)^2$ such that $\H_3 : x = 0$ and $\H_5 : y = 0$. Then $\H_7 : x + y = 0$ since $7 \equiv 3\cdot 5 \bmod 8$.
\begin{center}
    \includegraphics[scale=.2]{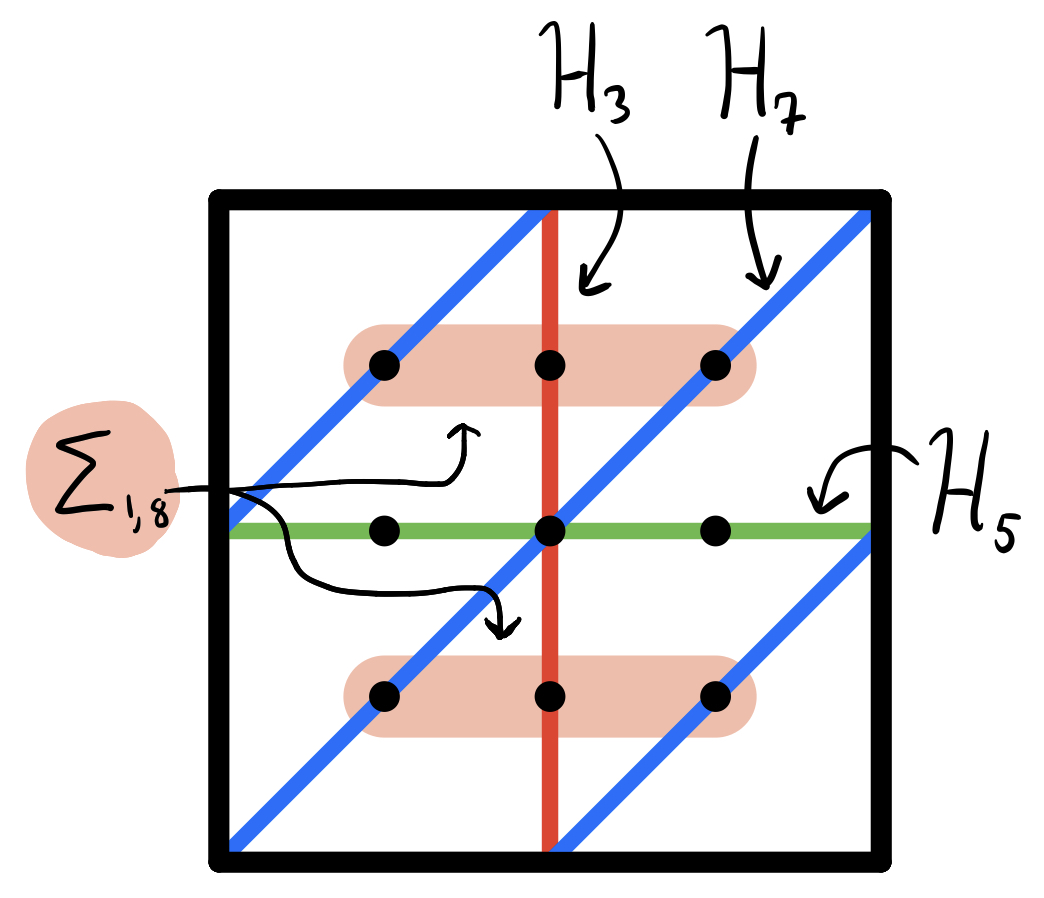}
\end{center}
Note that $\H_3 \cup \H_7$ does not contain the character $\chi \in \Sigma_{1,4}$ of conductor 4 which has additive coordinates $(1,0)$ in the picture above; this explains why $M_{21}(\zeta_4) \neq 0$.
\end{example}

\begin{example}
If $m = 6$, then $\h{\U}_6 = \{1, \chi\}$ contains only two characters. The non-trivial character $\chi$ has $c_\chi = 3$ and satisfies $\chi(2) = -1$. Thus the affine hyperplane $\H_2^a = \{\chi\}$ consists of the one non-trivial character. Therefore $\bigcup_{p\mid d}\H_p \cup \H_2^a$ covers $\h{\U}_6$ for any $d > 1$ coprime to $6$. Theorem \ref{thm intro precise vanish} implies that $M_{2d}(\zeta_6) = 0$ for all $d$ coprime to $6$. For example, this explains the $\Phi_6(x)$ factor in $M_{10}(x)$,
\[
    M_{10}(x) = g(x) \cdot \Phi_6\cdot \Phi_4\cdot \Phi_2\cdot \Phi_1\cdot x,
\]
where $g(x) = \tfrac{1}{10}(x^3 + x^2 - 1)$ is irreducible and not cyclotomic.
\end{example}

\begin{example}
Suppose we want to find some $d$ such that $\Phi_d(\zeta_8) = 1$. Since $m = 8$ is not divisible by 3 and $m' = 1$, we first look for $d$ such that $\H_{-1} \subseteq \bigcup_{p\mid 2d} \H_p$ in $\h{\U}_8 \cong \ZZ/(2)^2$. The hyperplane $\H_{-1}$ is one dimensional, hence may only be covered by itself. That is, $d$ must be divisible by a prime $p$ such that $p \equiv -1 \bmod 8$. The first example of such a $d$ is $d = 7$. However, $8$ does not divide $\varphi(7)$, whence $\Phi_7(\zeta_8) \neq 1$ by Theorem \ref{thm intro precise vanish}.2. If $d = 4991 = 7\cdot 23 \cdot 31$, then $8$ does divide $\varphi(4991)$ but $\varphi(4991)/8 = 495 \equiv 1 \bmod 2$ and
\[
    \sum_{a \mid 4991}\lfloor a/8\rfloor = 764 \equiv 0 \bmod 2,
\]
so again $\Phi_{4991}(\zeta_8) \neq 1$. If $d = 234577 = 7\cdot 23 \cdot 31 \cdot 47$, then we have that $8$ divides $\varphi(234577)$ and
\[
    \sum_{a\mid 234577}\lfloor a/8\rfloor = 36856 \equiv 0 \equiv 22770 = \frac{\varphi(234577)}{8} \bmod 2.
\]
Therefore Theorem \ref{thm intro precise vanish}.2 implies that $\Phi_{234577}(\zeta_8) = 1$.
\end{example}

\subsubsection{Necklace operators}
The connection between the necklace and shifted cyclotomic polynomials traces back to what we call the \emph{necklace operators} $\varphi_d$. Let $\NN^\circ$ denote the multiplicative semigroup of natural numbers, and let $\ZZ[\NN^\circ]$ be the integral semigroup ring comprised of all integral linear combinations of formal expressions $[m]$ with $m \in \NN$ subject only to the relations $[m][n] = [mn]$. The \emph{$d$th necklace operator} is defined by
\[
    \varphi_d := \sum_{e\mid d}\mu(e)[d/e] \in \ZZ[\NN^\circ].
\]
The polynomial ring $\QQ[x]$ carries a $\ZZ[\NN^\circ]$-module structure where $\alpha = \sum_m a_m [m] \in \ZZ[\NN^\circ]$ acts on $f(x) \in \QQ[x]$ by
\[
    \alpha f(x) := \sum_m a_m f(x^m).
\]
Similarly, the non-zero rational functions $\QQ(x)^\times$ have a multiplicative action of $\ZZ[\NN^\circ]$ defined on $g(x) \in \QQ(x)^\times$ by
\[
    g(x)^\alpha := \prod_m g(x^m)^{a_m}.
\]
With respect to these module structures we have the following expressions for necklace and cyclotomic polynomials in terms of the necklace operator,
\begin{align*}
    M_d(x) &= \frac{1}{d}\sum_{e\mid d}\mu(e)x^{d/e} = \frac{1}{d}\sum_{e\mid d}\mu(e)[d/e]x = \frac{\varphi_dx}{d}\\
    \Phi_d(x) &= \prod_{e\mid d} (x^{d/e} - 1)^{\mu(e)} = \prod_{e\mid d} (x - 1)^{\mu(e)[d/e]} = (x - 1)^{\varphi_d}.
\end{align*}
In Section \ref{sec neck op} we show how the abundance of pairs $(d,m)$ for which $M_d(\zeta_m) = 0$ or $\Phi_d(\zeta_m) = 1$ is ultimately a consequence of the elementary observation that the $d$th necklace operator has the following factorization in $\ZZ[\NN^\circ]$,
\begin{equation}
\label{eqn intro phi factored}
    \varphi_d = \prod_p [p^{m_p-1}]([p] - 1) = [d]\prod_{p \mid d}(1 - [p]^{-1}),
\end{equation}
where $d = \prod_p p^{m_p}$ is the prime factorization of $d$.

If $d$ is coprime to $m$, then $\varphi_d$ determines an element of the group ring $\ZZ[\s\U_m]$. Hence if $v$ is a vector in a linear $\U_m$-representation $V$, then $\varphi_d\s v \in V$. Our analysis of the identities $M_d(\zeta_m) = 0$ and $\Phi_d(\zeta_m) = 1$ hinges on the following result.

\begin{thm}
\label{thm intro abstract unlikely}
Let $d, m \geq 1$ be coprime integers and suppose $v \in V$ is an element of a $\QQ[\s\U_m]$-module. Let $\Sigma_v$ denote the set of Dirichlet characters that occur in the irreducible decomposition of the cyclic $\U_m$-representation generated by $v$. Then in $\CC \otimes V$ we have
\[
    \varphi_d\s v = \sum_{\chi \in \Sigma_v} \chi(d)\prod_{p\mid d}(1 - \overline{\chi(p)})v_\chi,
\]
where $v_\chi$ is the $\chi$-isotypic component of the vector $v$. Thus $\varphi_d\s v = 0$ if and only if $\Sigma_v \subseteq \bigcup_{p\mid d} \H_p$.
\end{thm}

\subsection{Related work}

As noted above, Necklace polynomials have many interpretations. Gauss \cite[Pg. 611]{Gauss} wrote down the necklace polynomials evaluated at a prime $p$ to count irreducible polynomials over $\FF_p$ of a prescribed degree and Sch\"onemann \cite[Sec. 48, Pp. 51-52]{Schon} later independently rediscovered this formula. This interpretation accounts for the appearance of necklace polynomials in the Euler product formula for the Hasse-Weil zeta function of the affine line over $\FF_q$,
\[
    \zeta_{\AA^1(\FF_q)}(t) = \frac{1}{1 - qt} = \prod_{d \geq 1}\left(\frac{1}{1 - t^d}\right)^{M_d(q)}.
\]
The name ``necklace polynomial'' comes from the combinatorial interpretation of $M_d(k)$ as counting the number of aperiodic necklaces of $d$ beads chosen from among $k$ colors, which Metropolis and Rota \cite[Pg. 95]{MR} attribute to the French colonel Moreau; the $M$ in the notation is presumably in his honor. Necklace polynomials also count Lyndon words \cite[Sec. 4.2]{BP} and the number of periodic orbits of a prescribed length for a generic polynomial of fixed degree \cite[Rmk. 4.3]{Silverman}. Metropolis and Rota \cite{MR} use necklace polynomials to construct a combinatorial model of the ring of big Witt vectors.

If $x = g$ is a natural number, then Witt \cite[Satz 3]{Witt} showed that $M_d(g)$ is the dimension of the degree $d$ homogeneous component of the free Lie algebra on $g$ generators. In this context the explicit expression for $M_d(x)$ as a divisor sum is sometimes called Witt's formula \cite[Pg. 1005]{BP}. Reutenaur \cite[Thm. 4.9, Thm. 5.1]{Reutenaur} gave a combinatorial proof of this result by constructing an explicit basis for the free Lie algebra from Lyndon words.

Let $\PConf_d(\RR^n)$ denote the space of labelled configurations of $d$ distinct points in $\RR^n$. The symmetric group $S_d$ acts naturally on this space by permuting labels and this action endows the cohomology $H^*(\PConf_d(\RR^n),\QQ)$ with the structure of an $S_d$-representation. The character values of these representations are determined by necklace polynomials. See Hyde \cite{hyde_poly}.

In \cite{hyde_ec}, we show that the values $M_d(\pm 1)$ of necklace polynomials at first and second order roots of unity may be interpreted as compactly supported Euler characteristics of spaces of degree $d$ irreducible polynomials over $\CC$ and $\RR$, respectively. In these cases the fundamental theorem of algebra gives a higher level explanation for why $M_d(\pm 1) = 0$ for nearly all $d$. It would be interesting to find a more conceptual interpretation of the vanishing of $M_d(\zeta_m)$ for $m > 2$, but we are unaware of one at this time.

The Euler characteristic interpretation of $M_d(\pm 1)$ found in \cite{hyde_ec} extends to the family $M_{d,n}(x)$ of \emph{higher necklace polynomials} introduced by the author in \cite{hyde_higher} to enumerate the irreducible polynomials over $\FF_q$ in $n$-variables. Theorem 1.5 in \cite{hyde_ec} shows that $M_{d,n}(\zeta_p) = 0$ for certain primes $p$ depending on $n$ and nearly all $d$. However, for $n > 1$, the qualitative behavior of these cyclotomic factors differs from those of $M_d(x)$ and $\Phi_d(x) - 1$, thus we expect the cyclotomic factors of $M_{d,n}(x)$ with $n > 1$ arise for a fundamentally different reason.

Despite the long history of necklace polynomials, the observation of their abundance of cyclotomic factors appears to be new.

The identity $\Phi_d(\zeta_m) = 1$ has received more attention. Note that if $\Phi_d(\zeta_m) = 1$, then
\begin{equation}
\label{eqn cyclo unit rel}
    1 = \Phi_d(\zeta_m) = \prod_{\gcd(j,d) = 1}(\zeta_m - \zeta_d^j).
\end{equation}
Algebraic integral units of the form $\zeta_m - \zeta_n$ are called \emph{cyclotomic units}. Thus \eqref{eqn cyclo unit rel} may be interpreted as a multiplicative relation between cyclotomic units. Such multiplicative relations are of interest in number theory and algebraic $K$-theory; they have been studied by Bass \cite{Bass}, Conrad \cite{conrad}, Ennola \cite{ennola}, Ramachandra \cite{ramachandra}, and others. This previous work focuses primarily on finding explicit relations that generate all of the relations amongst the cyclotomic units; our results provide a natural way of generating such relations through the construction of arrangements in $\h{\U}_m$ covering a prescribed set.

There is also some literature on classifying the vanishing integral linear combinations of roots of unity of which $M_d(\zeta_m) = 0$ and $\Phi_d(\zeta_m) - 1 = 0$ provide examples. See Christie, Dykema, Klep \cite{CDK} for a recent reference along with a survey of the previous work on this problem.

Kurshan and Odlyzko \cite{KO_applied, KO} made a detailed study of the unit part of $\Phi_d(\zeta_m)$ which included analyzing situations where $\Phi_d(\zeta_m) = 1$. Their work was motivated by problems related to the design of recursive linear digital filters. Our Theorem \ref{thm intro precise vanish}.2, characterizing solutions of $\Phi_d(\zeta_m) = 1$, is substantively equivalent to a result they proved in \cite{KO}. This can be seen most clearly in their discussion following Proposition 3.4 where they express the condition of $\H_{-1}$ being covered by hyperplanes as a disjunction of conditions on a character $\chi$ such that $\chi(-1) = 1$.

For their application, Kurshan and Odlyzko focus on analyzing the case $\Phi_{pm}(\zeta_m)$ with $p$ a prime not dividing $m$. This case is not covered by Theorem \ref{thm intro precise vanish}.2 since we assume that $m$ does not divide $def$.

The observation of the abundance and structure of cyclotomic factors of $\Phi_d(x) - 1$ and their connection to cyclotomic factors of $M_d(x)$ appears to be new.

The expression for $\varphi_d v$ given by Theorem \ref{thm key} generalizes a result of Bzd\k{e}ga, Herrera-Poyatos, Moree \cite[Thm. 1]{bzdega} which is the specialization to the case $v = \zeta_m - 1$ in the $\QQ$-linearization of $V = \QQ(\zeta_m)^\times$. They use this formula to explicitly evaluate $\Phi_d(\zeta_m)$ for small fixed values of $m$ as a function of $d$, which in turn they apply to give a new proof of a result of Vaughn on the heights of cyclotomic polynomials.

An earlier version of this manuscript appeared as Chapter 4 in the author's dissertation \cite{hyde_thesis}.

\subsection{Organization}
In Section \ref{sec neck op} we develop the theory of the necklace operators and prove Theorem \ref{thm intro abstract unlikely} as Theorem \ref{thm key}, from which we then deduce Theorem \ref{thm simple intro} as Corollary \ref{thm simple} and Theorem \ref{thm primewise intro} as Corollary \ref{cor primewise}. The first part of Theorems \ref{thm intro precise vanish} is proved in Section \ref{sec neck evals} as Theorem \ref{thm unlikely neck} and the second part is proved in Section \ref{subsec cyclo evals} as Theorem \ref{thm unlikely cyclo}.

\subsection{Acknowledgements}
We thank Andrew O'Desky for suggesting the connection with Dirichlet characters, for explaining an alternative proof of Proposition \ref{prop zeta support} via Gauss sums (see Remark \ref{remark gauss sum},) for many helpful, insightful conversations, and for extensive feedback on several drafts of this paper. 

We also thank Weiyan Chen, David Cox, Suki Dasher, Nir Gadish, Jeff Lagarias, Bob Lutz, and Phil Tosteson for their comments and feedback on earlier versions of the manuscript. The author is partially supported by the NSF MSPRF and the Jump Trading Mathlab Research Fund.

\section{Necklace operators}
\label{sec neck op}

We briefly review the representation theory of finite abelian groups---see Serre \cite{serre_rep} for more background. Given a finite (multiplicative) abelian group $\U$, let $\h{\U}$ denote the dual group or group of characters $\chi: \U \rightarrow \CC^\times$. The groups $\U$ and $\h{\U}$ are non-canonically isomorphic. Each character $\chi \in \h{\U}$ extends linearly to a ring homomorphism $\chi: \ZZ[\U] \rightarrow \CC$. If $\chi_i$ for $1 \leq i \leq n$ are the distinct characters of $\U$, then the map $\ZZ[\U] \rightarrow \CC^n$ given by
\[
    \alpha \in \ZZ[\U] \longmapsto (\chi_1(\alpha), \chi_2(\alpha), \ldots, \chi_n(\alpha)) \in \CC^n
\]
is an embedding of rings. Hence $\alpha \in \ZZ[\U]$ is zero if and only if $\chi(\alpha) = 0$ for all $\chi \in \h{\U}$. A \emph{hyperplane} $\H \subseteq \h{\U}$ is defined to be the (multiplicative) kernel of a character of $\h{\U}$. The group $\U$ is canonically isomorphic to the dual of $\h{\U}$. In particular, if $q \in \U$, then the hyperplane associated to $q$ is
\[
    \H_q := \ker(q) = \{\chi \in \h{\U} : \chi(q) = 1\}.
\]
If $q = 1$ is the identity, then $\H_1 = \h{\U}$ is the \emph{trivial hyperplane}. If $q \neq 1$, then $\H_q$ is a proper subgroup of $\h{\U}$.

\begin{remark}
\label{ex hyperplanes}
While we are primarily interested in multiplicative groups of units, the geometric terminology is best understood from an additive perspective. Suppose $\h{\U} \cong \FF_p^n$ is an $n$-dimensional vector space over a finite field $\FF_p$. If we choose some isomorphism of the $p$th roots of unity with the additive group of the field $\FF_p$, then a character $q: \h{\U} \rightarrow \CC^\times$ of $\h{\U}$ is equivalent under this isomorphism to an $\FF_p$-linear map $q: \FF_p^n \rightarrow \FF_p$. Thus there is a homogeneous linear form
\[
    h_q := \sum_{i=1}^n a_i x_i
\]
with $\FF_p$-coefficients such that the hyperplane $\H_q$ is precisely the set of solutions $h_q(x) = 0$ in $\FF_p^n$.
\end{remark}

Let $\U_m$ denote the group of units modulo $m$,
\[
    \U_m := (\ZZ/(m))^\times.
\]
The elements of $\h{\U}_m$ are called \emph{Dirichlet characters of modulus $m$}. If $n$ divides $m$, then the quotient map $\U_m \rightarrow \U_n$ induces an injective map $\h{\U}_n \rightarrow \h{\U}_m$. Identifying $\h{\U}_n$ with its image under this map we say $\h{\U}_n \subseteq \h{\U}_m$. If a character $\chi \in \h{\U}_m$ belongs to the subset $\h{\U}_n$, then we say $\chi$ has \emph{modulus} $n$. If $\chi$ has modulus $n$, then the values $\chi(k)$ depend only on $k$ modulo $n$. Note that if $\chi$ has modulus $n$, it also has modulus $m$ for all multiples $m$ of $n$. The smallest $n$ for which $\chi \in \h{\U}_m$ has modulus $n$ is called the \emph{conductor} of $\chi$ and denoted $c_\chi$.

\begin{caution}
\label{caution}
A common convention in number theory is to distinguish a character $\chi \in \h{\U}_n$ from the character it naturally induces in $\h{\U}_m$ when $n \mid m$. In particular, the convention is to set $\chi(d) = 0$ for all non-trivial $\chi \in \h{\U}_m$ when $d$ is not coprime to $m$. Since we are identifying $\h{\U}_n$ with a subset of $\h{\U}_m$ whenever $n$ divides $m$, we use a slight natural variation on this convention: If $d \in \ZZ$ and $\chi \in \h{\U}_m$ has conductor $n$, then we set $\chi(d) = 0$ if $d$ is not coprime to $n$ and otherwise set $\chi(d)$ to the well-defined, nonzero value of $\chi$ on the residue class of $d$ modulo $n$. This gives each character $\chi$ a consistent value independent of which group $\h{\U}_m$ it is considered to be an element of. Our convention will prove to be a useful simplification throughout this paper.
\end{caution}

\begin{example}
\label{ex caution}
If $m = 10$, then $\h{\U}_{10} = \h{\U}_5$. If $\chi \in \h{\U}_{10}$, then the common convention is to say that $\chi(2) = 0$ since $2$ divides $10$. However, $\chi$ has conductor 5 and as an element of $\h{\U}_5$ it has a well-defined non-zero value at $2$ which we take to be the value of $\chi(2)$.
\end{example}

If $R$ is a semiring, we let $R^\circ$ denote the multiplicative semigroup of $R$. Let $\ZZ[\NN^\circ]$ denote the ring generated by expressions $[m]$ with $m \in \NN$ subject to the relations $[m][n]= [mn]$. We define the \emph{$d$th necklace operator} for $d\geq 1$ to be the element $\varphi_d \in \ZZ[\NN^\circ]$ defined by
\[
    \varphi_d := \sum_{e\mid d} \mu(e)[d/e].
\]

\begin{remark}
The map $[n] \mapsto n$ determines a ring homomorphism $\ZZ[\NN^\circ] \rightarrow \ZZ$ such that
\[
    \varphi_d \longmapsto \sum_{e\mid d} \mu(e)(d/e) = \varphi(d),
\]
where $\varphi(d)$ is the \emph{Euler totient function}, hence our choice of notation.
\end{remark}

Necklace polynomials and cyclotomic polynomials are connected through the necklace operator. Recall from the introduction that with respect to the natural additive and multiplicative actions of $\ZZ[\NN^\circ]$ on $\QQ[x]$ and $\QQ(x)^\times$, respectively, we have
\[
    M_d(x) = \frac{\varphi_d x}{d}, \hspace{.45in} \Phi_d(x) = (x - 1)^{\varphi_d}.
\]

The map $[n] \mapsto [n \bmod m]$ induces a ring homomorphism $\ZZ[\NN^\circ] \rightarrow \ZZ[\ZZ/(m)^\circ]$. If $d$ is coprime to $m$, then the image of $\varphi_d$ under this map belongs to the subring $\ZZ[\s\U_m]$.
The image of $\varphi_d$ in $\ZZ[\s\U_m]$ factors as
\begin{equation}
\label{eqn phi factored}
    \varphi_d = [d]\prod_{p\mid d}(1 - [p]^{-1}).
\end{equation}
The factorization \eqref{eqn phi factored} is equivalent to families of functional identities satisfied by $M_d(x)$ and $\Phi_d(x)$: If $p$ is a prime and $d\geq 1$, then
\[
    M_{dp}(x) = \begin{cases} \tfrac{1}{p}(M_d(x^p) - M_d(x)) & p \nmid d\\ \tfrac{1}{p}M_d(x^p) & p \mid d \end{cases}, \hspace{.45in}
    \Phi_{dp}(x) = \begin{cases} \Phi_d(x^p)/\Phi_d(x) & p \nmid d\\ \Phi_d(x^p) & p \mid d \end{cases}.
\]
The identities for necklace polynomials were observed and given combinatorial interpretations by Metropolis, Rota \cite{MR}; the identities for cyclotomic polynomials are well-known.

\begin{remark}
\label{rem sqfree}
Let $d_0$ be the product of all distinct primes dividing $d$. Thus \eqref{eqn phi factored} implies that $\varphi_d = [d/d_0]\varphi_{d_0}$, hence
\begin{equation}
\label{eqn core reduction}
    dM_d(x) = d_0M_{d_0}(x^{d/d_0}), \hspace{.45in} \Phi_d(x) = \Phi_{d_0}(x^{d/d_0}).
\end{equation}
We use \eqref{eqn core reduction} to reduce the analysis of $M_d(\zeta_m)$ and $\Phi_d(\zeta_m)$ to the case where $d$ is squarefree.
\end{remark}

Let $\CC[\s \U_m]$ denote the group algebra of $\U_m$ over $\CC$. If $\chi \in \h{\U}_m$ is a character, let $e_\chi \in \CC[\s\U_m]$ denote the corresponding idempotent,
\[
    e_\chi := \frac{1}{\varphi(m)}\sum_{q \in \s\U_m}\overline{\chi(q)}[q].
\]
We write $v_\chi := e_\chi v$ for the projection of a vector $v \in V$ onto the $\chi$-isotypic component of $V$. Then
\[
    v = \sum_{\chi \in \,\h{\U}_m} v_\chi.
\]
The \emph{support of $v$} is the set $\Sigma_v \subseteq \h{\U}_m$ of characters $\chi$ such that $v_\chi \neq 0$. In particular, $v = 0$ if and only if $\Sigma_v = \emptyset$.

\begin{thm}
\label{thm key}
Let $d, m \geq 1$ be coprime integers, let $V$ be a $\QQ[\s\U_m]$-module, and let $v \in V$ be an element with support $\Sigma_v$. Then $\varphi_dv$ has the following expression in $\CC \otimes V$,
\[
    \varphi_d\s v = \sum_{\chi \in \Sigma_v} \chi(d)\prod_{p\mid d}(1 - \overline{\chi(p)})v_\chi.
\]
Thus $\varphi_d\s v = 0$ if and only if $\Sigma_v \subseteq \bigcup_{p\mid d} \H_p$. In particular, $\varphi_d = 0$ in $\QQ[\s\U_m]$ if and only if \,$\h{\U}_m \subseteq \bigcup_{p\mid d}\H_p$.
\end{thm}

\begin{proof}
If $\alpha \in \CC[\s\U_m]$ and $e_\chi$ is the idempotent associated to a character $\chi$, then $\alpha e_\chi = \chi(\alpha) e_\chi$. Hence
\begin{equation}
\label{eqn phi image id}
    \varphi_d\s v = \sum_{\chi \in \Sigma_v}\varphi_d e_\chi v
    = \sum_{\chi \in \Sigma_v} \chi(\varphi_d) v_\chi
    = \sum_{\chi \in \Sigma_v} \chi(d)\prod_{p\mid d}(1 - \overline{\chi(p)})v_\chi
\end{equation}
where the final equality follows from \eqref{eqn phi factored}. The factor $1 - \overline{\chi(p)}$ vanishes precisely when $\chi \in \H_p$, thus the support of $\varphi_d\s v$ is $\Sigma_v \setminus \bigcup_{p\mid d}\H_p$. Therefore $\varphi_dv = 0$ if and only if $\Sigma_v \setminus \bigcup_{p\mid d} \H_p = \emptyset$, which is to say $\Sigma_v \subseteq \bigcup_{p\mid d} \H_p$. Since $\QQ[\s\U_m]$ is cyclic as a module over itself generated by $1$, $\Sigma_1 = \h{\U}_m$ and it follows that $\varphi_d = \varphi_d 1 = 0$ if and only if $\h{\U}_m \subseteq \bigcup_{p\mid d}\H_p$.
\end{proof}

Theorem \ref{thm key} gives us the following simple sufficient condition for both of the identities $M_d(\zeta_m) = 0$ and $\Phi_d(\zeta_m) = 1$ to hold simultaneously.

\begin{cor}
\label{thm simple}
Let $d, m > 1$ be coprime integers. If $\h{\U}_m \subseteq \bigcup_{p\mid d} \H_p$,
then $x^m - 1$ divides $M_d(x)$ and $\frac{x^m-1}{x-1}$ divides $\Phi_d(x) - 1$.
\end{cor}

\begin{proof}
If $\h{\U}_m \subseteq \bigcup_{p\mid d}\H_p$, then $\varphi_d = 0$ by Theorem \ref{thm key}. Thus,
\[
    M_d(\zeta_m^k) = \frac{\varphi_d \zeta_m^k}{d} = 0\hspace{.45in} \Phi_d(\zeta_m^k) = (\zeta_m^k - 1)^{\varphi_d} = 1.
\]
The first identity holds for all $k\geq 0$, but in the second identity we need $k\not\equiv 0 \bmod m$ in order for $\zeta_m^k - 1 \in \QQ(\zeta_m)^\times$. Therefore $x^m - 1$ divides $M_d(x)$ and $\frac{x^m - 1}{x - 1}$ divides $\Phi_d(x) - 1$.
\end{proof}

\begin{example}
\label{ex lines}
Let $m = 65$ and let $d = 9372603371 = 47 \cdot 73 \cdot 79 \cdot 151 \cdot 229$. The group $\h{\U}_{65}$ decomposes as $\h{\U}_{65} \cong \ZZ/(4)^2 \times \ZZ/(3)$, hence each hyperplane $\H_p$ factors as $\H_p \cong \H_p^{(4)} \times \H_p^{(3)}$ with $\H_p^{(4)} \subseteq \ZZ/(4)^2$ and $\H_p^{(3)} \subseteq \ZZ/(3)$. In this case, each of the hyperplanes $\H_p$ with $p \mid d$ is trivial in the 3-component $\H_p^{(3)} = \ZZ/(3)$. Thus we can visualize the hyperplanes $\H_p$ via their 4-component $\H_p^{(4)}$ as lines in the ``plane'' $\ZZ/(4)^2$. Each of the five primes dividing $d$ corresponds to a different colored line in the diagram below with respect to the choice of coordinates $x = \rho(47)$ and $y = \rho(151)$. Since the five lines $\H_p$ with $p\mid d$ cover all of $\h{\U}_{65}$, Corollary \ref{thm simple} implies that $M_d(\zeta_{65}^k) = 0$ for all $k\geq 0$ and $\Phi_d(\zeta_{65}^k) = 1$ for all $k\not\equiv 0 \bmod 65$.

\begin{center}
    \includegraphics[scale=.2]{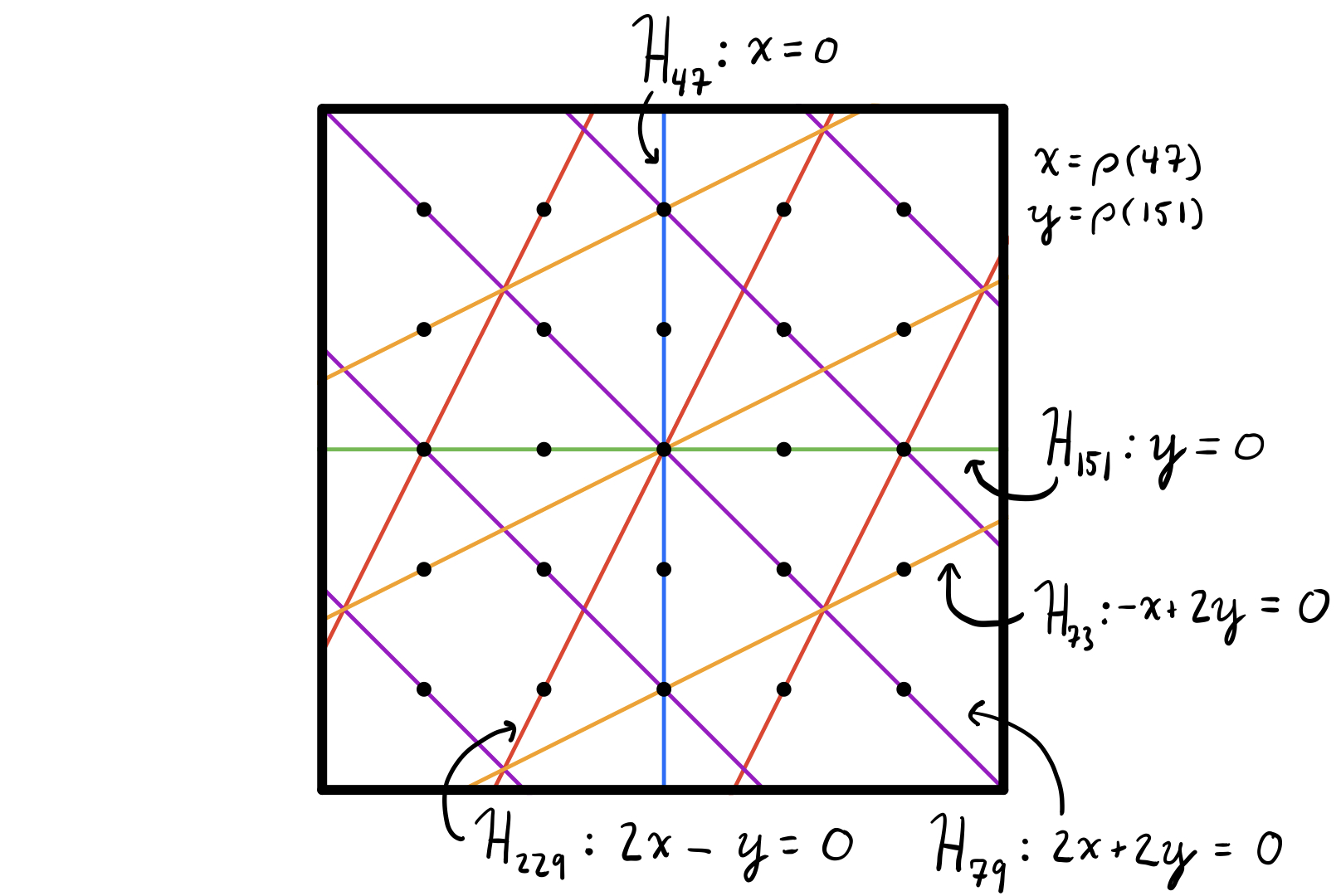}
\end{center}

By drawing other arrangements of lines covering $\ZZ/(4)^2$ and then finding primes in the corresponding congruence classes modulo $65$ (which exist by Dirichlet's theorem on primes in arithmetic progressions) we can construct several other nontrivial examples of $d$ for which $M_d(\zeta_{65}) = 0$ and $\Phi_d(\zeta_{65}) = 1$.
\begin{center}
    \includegraphics[scale=.2]{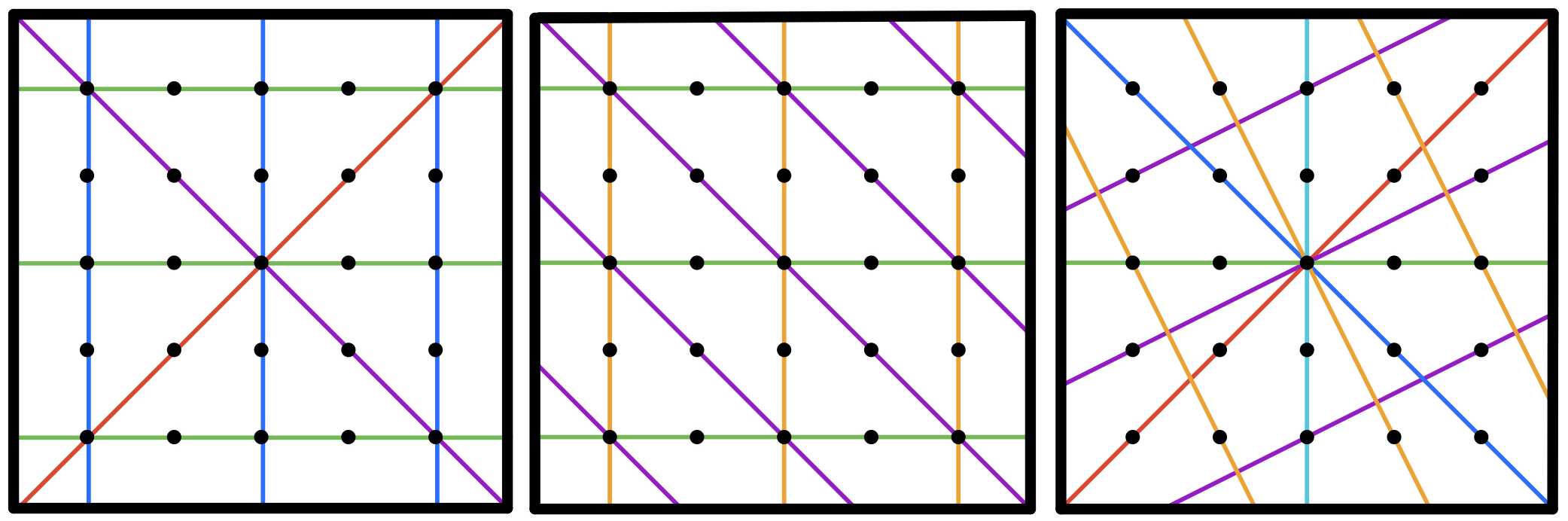}
\end{center}
Example values of $d$ for each of these arrangements are, respectively,
\begin{align*}
    d_1 &= 157\cdot 181 \cdot 337 \cdot 389\\
    d_2 &= 79 \cdot 181 \cdot 389\\
    d_3 &= 47 \cdot 109 \cdot 151 \cdot 157 \cdot 317 \cdot 337.\qedhere
\end{align*}
\end{example}

The following corollary of Theorem \ref{thm key} proves Theorem \ref{thm primewise intro} from the introduction.

\begin{cor}
\label{cor primewise}
Let $d, e, m \geq 1$ be integers and let $v \in V$ be an element of a $\QQ[\s\U_m]$-module $V$.
\begin{enumerate}
    \item If $\varphi_dv = 0$ and $e$ is coprime to $m$, then $\varphi_{de}v = 0$.
    \item If $d$ and $e$ are coprime to $m$ and the following sets are equal
    \[
        \{p \bmod m : p \mid d \text{ is prime}\} = \{q \bmod m : q \mid e \text{ is prime}\},
    \]
    then $\varphi_d v = 0$ if and only if $\varphi_ev = 0$.
\end{enumerate}
In particular, if $V = \QQ \otimes \QQ(\zeta_m)$ and $v = \zeta_m$; or if $m > 1$,  $V = \QQ(\zeta_m)^\times$, and $v = \zeta_m - 1$, then (1) and (2) hold with $\varphi_dv = M_d(\zeta_m)$ and $\varphi_dv = \Phi_d(\zeta_m)$, respectively.
\end{cor}

\begin{proof}
(1) The product formula \eqref{eqn phi factored} for the necklace operator implies that $\varphi_d$ divides $\varphi_{de}$ in $\ZZ[\NN^\circ]$ and the assumption that $e$ is coprime to $m$ implies that $\varphi_{de}/\varphi_d \in \ZZ[\s\U_m]$. Thus,
\[
    \varphi_{de}v = (\varphi_{de}/\varphi_d)(\varphi_dv) = (\varphi_{de}/\varphi_d)0 = 0.
\]

(2) Theorem \ref{thm key} implies that $\varphi_dv = 0$ if and only if  $\Sigma_v \subseteq \bigcup_{p\mid d}\H_p$, where $\Sigma_v$ is the support of $v$. The hyperplane $\H_p \subseteq \h{\U}_m$ depends only on the residue class $p \bmod m$. In other words,
\begin{equation}
\label{eqn sets}
    \{p \bmod m : p \mid d \text{ is prime}\} = \{q \bmod m : q \mid e \text{ is prime}\}
\end{equation}
is equivalent to $\bigcup_{p\mid d} \H_p = \bigcup_{q\mid e}\H_q$. Therefore, \eqref{eqn sets} implies  $\varphi_dv = 0$ if and only if $\varphi_e v = 0$.
\end{proof}


\section{Cyclotomic factors of \texorpdfstring{$M_d(x)$}{Md(x)}}
\label{sec neck evals}

In this section we characterize those pairs $(d,m)$ for which $M_d(\zeta_m) = 0$ in terms of an explicit set of Dirichlet characters being covered by an arrangement of hyperplanes. If $d$ is coprime to $m$, this reduces to determining the support of $\zeta_m \in \QQ(\zeta_m)$ by Theorem \ref{thm key}. When $d$ and $m$ are not coprime, the situation becomes more complicated and the relevant support depends in a subtle way on the common factors of $d$ and $m$.

If $m$ is a positive integer, then the \emph{squarefree part} of $m$, denoted $m'$, is the product of all primes that divide $m$ exactly once. We say a character $\chi \in \h{\U}_m$ is \emph{supportive} if the conductor of $\chi$ is divisible by $m/m'$. Equivalently, $\chi$ is supportive if and only if
\begin{equation}
\label{eqn local}
    v_p(c_\chi) = v_p(m) \text{ for all primes $p$ such that $v_p(m) \geq 2$.}
\end{equation}
Let $\h{\U}_m^* \subseteq \h{\U}_m$ denote the subset of all supportive characters.

\begin{prop}
\label{prop zeta support}
Let $m \geq 1$ be an integer and let $\chi \in \h{\U}_m$ be a character of modulus $m$. Then $\chi$ is in the support of $\zeta_m \in \QQ(\zeta_m)$ if and only if $\chi$ is supportive.
\end{prop}

\begin{proof}
If $m = \prod_p p^{m_p}$ is the prime factorization of $m$ and $\chi \in \h{\U}_m$ is a character, then by the Chinese Remainder Theorem there are factorizations
\[
    \zeta_m = \prod_p \zeta_{p^{m_p}}^{a_p}, \hspace{.25in} \chi = \prod_p \chi_{p^{m_p}},
\]
where $a_p \in \U_{p^{m_p}}$ is some unit and $\chi_{p^{m_p}} \in \h{\U}_{p^{m_p}}$ is some character of modulus $p^{m_p}$. The factorization of $\chi$ induces a factorization of idempotents $e_\chi = \prod_p e_{\chi_{p^{m_p}}}$ such that
\begin{equation}
\label{eqn idem prod decomp}
    e_\chi \zeta_m = \prod_p e_{\chi_{p^{m_p}}}\zeta_{p^{m_p}}^{a_p}.
\end{equation}
Let $\Sigma_m$ denote the support of $\zeta_m$ in $\QQ(\zeta_m)$. Since $\Sigma_m$ depends only on the cyclic $\QQ[\s\U_m]$-module generated by $\zeta_m$, it follows that $\Sigma_m$ is the support of $\zeta_m^a$ for all $a \in \U_m$. Then \eqref{eqn idem prod decomp} implies that $\chi \in \Sigma_m$ if and only if $\chi_{p^{m_p}} \in \Sigma_{p^{m_p}}$ for all primes $p$. The conductor of $\chi$ is the product of the conductors of $\chi_{p^{m_p}}$, hence by the definition of supportive characters, $\chi \in \h{\U}_m^*$ if and only if $\chi_{p^{m_p}} \in \h{\U}_{p^{m_p}}^*$. Thus to prove our claim it suffices to show that $\Sigma_{p^k} = \h{\U}_{p^k}^*$ for all primes $p$ and all $k\geq 1$. Note that $\h{\U}_{p^k}^* \subseteq \h{\U}_{p^k}$ consists of all the primitive characters if $k > 1$ and all characters if $k = 1$.

If $k = 1$, then the identity
\[
    \sum_{q \in \U_p} \zeta_p^q = -1
\]
implies that $\QQ(\zeta_p)$ is the cyclic $\QQ[\s\U_p]$-module generated by $\zeta_p$. Therefore $\Sigma_p = \h{\U}_p = \h{\U}_p^*$.

If $k > 1$, then $\{1, \zeta_{p^k}, \zeta_{p^k}^2,\ldots, \zeta_{p^k}^{p-1}\}$ forms a $\QQ(\zeta_{p^{k-1}})$-basis for $\QQ(\zeta_{p^k})$. Thus $\QQ(\zeta_{p^k})$ decomposes as the direct sum of the following two $\QQ[\s\U_{p^k}]$-submodules,
\begin{equation}
\label{eqn p decomp}
    \QQ(\zeta_{p^k}) = \QQ(\zeta_{p^{k-1}}) \oplus \sum_{a=1}^{p-1} \QQ(\zeta_{p^{k-1}})\zeta_{p^k}^a =: U \oplus V.
\end{equation}
We claim that $V$ is the cyclic $\QQ[\s\U_{p^k}]$-module generated by $\zeta_{p^k}$.
If $q \in \U_{p^k}$, then $q \equiv a + pb \bmod p^k$ for some $1\leq a \leq p - 1$ and some integer $b$. Thus $\zeta_{p^k}^q = \zeta_{p^{k-1}}^b \zeta_{p^k}^a \in \QQ(\zeta_{p^{k-1}})\zeta_{p^k}^a$, and elements of this form span $V$ by construction. The normal basis theorem implies that $\QQ(\zeta_m) \cong \QQ[\s\U_m]$ as $\QQ[\s\U_m]$-modules for any $m$. Recall that $\h{\U}_{p^{k-1}}$ is identified with its natural image in $\h{\U}_{p^k}$ (see Caution \ref{caution}.) Therefore, taking supports in \eqref{eqn p decomp} gives us
\[
    \h{\U}_{p^k} = \h{\U}_{p^{k-1}} \sqcup \Sigma_{p^k}.
\]
Therefore $\Sigma_{p^k} = \h{\U}_{p^k} \setminus \h{\U}_{p^{k-1}} = \h{\U}_{p^k}^*$ consists of the primitive characters of modulus $p^k$.
\end{proof}

\begin{remark}
\label{remark gauss sum}
If $\chi \in \h{\U}_m$ is a non-trivial Dirichlet character of modulus $m$, then the \emph{Gauss sum} of $\chi$ is
\[
    G(\chi) := \sum_{q\in \s\U_m}\chi(q)\zeta_m^q.
\]
Gauss sums are scalar multiples of the isotypic components of $\zeta_m \in \QQ(\zeta_m)$. In particular,
\[
    G(\chi^{-1}) = \sum_{q\in \s\U_m} \overline{\chi(q)}\zeta_m^q = \varphi(m)e_{\chi}\zeta_m = \varphi(m) (\zeta_{m })_{\chi}.
\]
Thus the support of $\zeta_m$ may be interpreted as the set of all characters $\chi$ such that $G(\chi^{-1}) \neq 0$. Since $\h{\U}_m^*$ is closed under taking inverses, Proposition \ref{prop zeta support} is equivalent to the assertion
\[
    \h{\U}_m^* = \{\chi \in \h{\U}_m : G(\chi) \neq 0\}.
\]
This characterization of non-vanishing Gauss sums, and hence of the support of $\zeta_m$, may also be deduced from the classical theory of Gauss sums. In particular, it follows from Theorems 9.7 and 9.10 in Montgomery, Vaughn \cite{Montgomery}. We thank Andrew O'Desky for bringing this to our attention.
\end{remark}

We now turn to the main result of this section.

\begin{thm}
\label{thm unlikely neck}
Let $d, e, f, m \geq 1$ be integers, let $m'$ be the squarefree part of $m$, and let $\H_2^a \subseteq \h{\U}_m$ be the affine hyperplane $\H_2^a := \{\chi \in \h{\U}_m : \chi(2) = -1\}$. Suppose that
\begin{multicols}{2}
\begin{enumerate}[label = (\roman*)]
    \item $def$ is squarefree,
    \item $d$ is coprime to $m$,
    \item $e$ divides $m'$,
    \item $f$ divides $m/m'$.
\end{enumerate}
\end{multicols}
Let $\Sigma_{f,m} \subseteq \h{\U}_m$ be the set of all characters $\chi$ such that
\begin{enumerate}
    \item $v_p(c_\chi) = v_p(m)$ if $v_p(m) \geq 2$ and $v_p(f) = 0$, and
    \item $v_p(c_\chi) \geq v_p(m) - 1$ if $v_p(m) > 2$ and $v_p(f) = 1$.
\end{enumerate}
Then $M_{def}(\zeta_m) = 0$ if and only if
\[
     \Sigma_{f,m}\subseteq \begin{cases} \bigcup_{p\mid d} \H_p & \text{if $e$ is odd,} \\ \bigcup_{p\mid d} \H_p \cup \H_2^a & \text{if $e$ is even.}\end{cases}
\]
\end{thm}

\begin{remark}
Recall that by Remark \ref{rem sqfree} we lose no generality in assuming that $def$ is squarefree.
\end{remark}

\begin{proof}
Since we assume $d, e, f$ are pairwise coprime we may express $M_{def}(\zeta_m)$ as
\[
    M_{def}(\zeta_m) = \frac{1}{def}\varphi_d(\varphi_e\varphi_f\zeta_m).
\]
Our strategy is to determine the support of $\varphi_e\varphi_f\zeta_m$ and then apply Theorem \ref{thm key} with $v = \varphi_e\varphi_f\zeta_m$. Note that Theorem \ref{thm key} does not immediately apply with $v = \zeta_m$ because $def$ is not coprime to $m$.

Observe that
\[
    \varphi_f\s\zeta_m = \sum_{b\,\mid\, f} \mu(f/b) \zeta_m^{b} = \sum_{b\,\mid\, f} \mu(f/b)\zeta_{m/b},
\]
where $\zeta_{m/b} := \zeta_m^{b}$ is a primitive $m/b$th root of unity and $\mu(f/b) \neq 0$ since $f$ is squarefree. Proposition \ref{prop zeta support} implies that the support of $\zeta_{m/b}$ is $\h{\U}_{m/b}^*$. If $b, b'$ are distinct divisors of $f$, then by the definition of $f$ there is some prime $p$ such that $v_p(m) \geq 2$ and, say, $1 = v_p(b) > v_p(b') = 0$. Thus if $c$ and $c'$ are the conductors of characters in $\h{\U}_{m/b}^*$ and $\h{\U}_{m/b'}^*$, respectively, then
\[
    v_p(c) \leq v_p(m/b) < v_p(m) = v_p(c'),
\]
where the last equality follows from Proposition \ref{prop zeta support}. In particular, $\h{\U}_{m/b}^*$ and $\h{\U}_{m/b'}^*$ are disjoint. Therefore the support of $\varphi_f\s\zeta_m$ is
\[
    \Sigma_{\varphi_f\zeta_m} = \bigcup_{b\,\mid\, f}\h{\U}_{m/b}^*
\]
Let $\Sigma_{f,m} \subseteq \h{\U}_m$ be the set of characters defined in the statement of Theorem \ref{thm unlikely neck}. We claim that
\begin{equation}
\label{eqn set simplify}
    \Sigma_{f,m} = \bigcup_{b\mid f}\h{\U}_{m/b}^* = \Sigma_{\varphi_f\zeta_m}.
\end{equation}

Suppose that $\chi \in \h{\U}_{m/b}^*$ for some $b \mid f$. Then \eqref{eqn local} implies that $v_p(c_\chi) = v_p(m/b)$ whenever $v_p(m/b) \geq 2$. Since $b$ is squarefree, there are two cases: if $v_p(b) = 0$, then $v_p(m) = v_p(m/b) \geq 2$ and $v_p(c_\chi) = v_p(m)$; and if $v_p(b) = 1$, then $v_p(m) > 2$ and $v_p(c_\chi) = v_p(m) - 1$. Hence $\chi \in \Sigma_{f,m}$ and thus $\bigcup_{b\mid f}\h{\U}_{m/b}^* \subseteq \Sigma_{f,m}$.

For the reverse inclusion, suppose that $\chi \in \Sigma_{f,m}$. Let $b$ be the product of all primes $p \mid f$ such that $v_p(c_\chi) < v_p(m)$. Then $b$ is a divisor of $f$ and $c_\chi$ divides $m/b$. If $p$ is a prime such that $v_p(b) = 0$, then $v_p(c_\chi) = v_p(m)$ by construction. If $p$ is a prime such that $v_p(b) = 1$ and $v_p(m/b) > 1$, then $v_p(m) > 2$ and $v_p(c_\chi) \geq v_p(m) - 1$ by the definition of $\Sigma_{f,m}$. On the other hand, $v_p(c_\chi) < v_p(m)$ since $p$ divides $b$, hence $v_p(c_\chi) = v_p(m) - 1$. In either case we have $v_p(c_\chi) = v_p(m/b)$ when $v_p(m/b) \geq 2$, which is equivalent to $\chi \in \h{\U}_{m/b}^*$. Therefore $\Sigma_{f,m} \subseteq \bigcup_{b\mid f} \h{\U}_{m/b}^*$, which finishes the proof of \eqref{eqn set simplify}.

Now suppose $p$ is a prime dividing $e$, so that $v_p(m) = 1$. Since $m/p$ is coprime to $p$ by assumption, we may write $\zeta_m = \zeta_{m/p}^a\zeta_p^b$ for some $a\in \U_{m/p}$ and $b\in \U_p$. Recall that
\[
    1 = -\sum_{k=1}^{p-1}\zeta_p^{bk}.
\]
For $1 \leq k \leq p-1$, let $c(p,k) \in \U_m$ be the unique unit such that
\begin{align*}
    c(p,k) &\equiv p \bmod m/p\\
    c(p,k) &\equiv k \bmod p.
\end{align*}
Then $\zeta_m^{c(p,k)} = \zeta_{m/p}^{ap}\zeta_p^{bk}$. Hence
\[
    \varphi_p\s\zeta_m = \zeta_m^{p} - \zeta_m
    = \zeta_{m/p}^{ap} - \zeta_m
    = - \Big(\zeta_m  + \sum_{k=1}^{p-1}\zeta_{m/p}^{ap}\zeta_p^{bk}\Big)
    = -\Big(1 + \sum_{k=1}^{p-1}[c(p,k)]\Big)\zeta_m 
    =: -\alpha_p\s\zeta_m.
\]

Recall that any $\chi \in \h{\U}_m$ can be factored as $\chi = \chi_{m/p}\chi_p$ with $\chi_n \in \h{\U}_n$ (see the proof of Proposition \ref{prop zeta support}.) Then
\[
    \chi(c(p,k)) = \chi_{m/p}(p)\chi_p(k).
\]
Hence
\[
    \chi(\alpha_p) = 1 + \chi_{m/p}(p)\sum_{k=1}^{p-1}\chi_p(k).
\]
The orthogonality relations for characters imply that
\[
    \sum_{k=1}^{p-1}\chi_p(k) = \begin{cases} p - 1 & \chi_p = 1,\\ 0 & \chi_p \neq 1.\end{cases}
\]
If $\chi_p \neq 1$, then $\chi(\alpha_p) = 1 \neq 0$. If $\chi_p = 1$, then $\chi(\alpha_p) = 0$ is equivalent to
\[
    0 = \chi(\alpha_p) = 1 + \chi_{m/p}(p)(p - 1) \Longrightarrow \chi_{m/p}(p) = \frac{1}{1 - p}.
\]
Since $\chi_{m/p}(p)$ is a root of unity, it must be the case that $p = 2$ and $\chi_{m/p}(2) = -1$. In other words, $\chi(\alpha_p) = 0$ if and only if $\chi \in \H_2^a$. Thus $\Sigma_{\varphi_p\zeta_m} = \Sigma_{\zeta_m} = \h{\U}_m^*$ for each odd prime $p\mid e$ and $\Sigma_{\varphi_2\zeta_m} = \h{\U}_m^* \setminus \H_2^a$. Since $e$ is squarefree, $\varphi_e = \prod_{p\mid e} \varphi_p$, hence
\[
    \Sigma_{\varphi_e\varphi_f\zeta_m} = \begin{cases} \Sigma_{f,m} & \text{if $e$ is odd,}\\ \Sigma_{f,m} \setminus \H_2^a & \text{if $e$ is even.}\end{cases}
\]
Thus our conclusion follows from Theorem \ref{thm key}.
\end{proof}

\begin{cor}
\label{cor m squarefree}
Let $d, e, m\geq 1$ be as in Theorem \ref{thm unlikely neck}. Suppose that $m$ is squarefree and $e$ is odd. If $M_{de}(\zeta_m) = 0$, then $M_{de}(\zeta_m^k) = 0$ for all $k\geq 0$. In other words, $\Phi_m(x)$ divides $M_{de}(x)$ if and only if $x^m - 1$ divides $M_{de}(x)$.
\end{cor}

\begin{proof}
If $m$ is squarefree, then $\h{\U}_m^* = \h{\U}_m$ and $f = 1$ in the notation of Theorem \ref{thm unlikely neck}. Since $e$ is odd, Theorem \ref{thm unlikely neck} implies that $M_{de}(\zeta_m) = 0$ if and only if $\bigcup_{p\mid d}\H_p$ covers $\h{\U}_m$, and this is equivalent to $\varphi_d = 0$ by Theorem \ref{thm key}. Thus for all $k\geq 0$,
\[
    M_{de}(\zeta_m^k) = \frac{\varphi_d}{d}M_e(\zeta_m^k) = 0.\qedhere
\]
\end{proof}

\begin{example}
In the introduction we observed that
\[
    M_{253}(x) = f(x)\cdot \Phi_{24}\cdot \Phi_{22}\cdot \Phi_{11}\cdot \Phi_{10}\cdot \Phi_8\cdot \Phi_5\cdot \Phi_2\cdot \Phi_1\cdot x
\]
for some non-cyclotomic irreducible polynomial $f(x) \in \QQ[x]$. Notice that for each squarefree $m$ such that $\Phi_m(x)$ divides $M_{253}(x)$ we have that $x^m - 1$ divides $M_{253}(x)$, but this property fails to hold for the non-squarefree $m = 24$; this reflects Corollary \ref{cor m squarefree}.

To see the necessity of the condition that $e$ is odd in Corollary \ref{cor m squarefree} consider the factorization
\[
    M_{10}(x) = g(x)\cdot \Phi_6\cdot \Phi_4\cdot \Phi_2\cdot \Phi_1\cdot x
\]
for some non-cyclotomic irreducible polynomial $g(x) \in \QQ[x]$. If $m = 6$, then $d = 5$ and $e = 2$; we see that $\Phi_6(x)$ divides $M_{10}(x)$, but $\Phi_3(x)$ does not.

In Example \ref{ex m = 8} we showed that $M_{21}(\zeta_8) = 0$ but $M_{21}(\zeta_8^2) \neq 0$, which shows the necessity of the $m$ squarefree assumption.
\end{example}

\begin{example}
Let $f = 3$ and let $m = 1026 = 2\cdot 3^3 \cdot 19$. Then the set $\Sigma_{f,m}$ of characters defined in Theorem \ref{thm unlikely neck} simplifies in this case to
\[
    \Sigma_{3,1026} = \{\chi \in \h{\U}_{1026} : v_3(c_\chi) \geq 2\}.
\]
If $\chi \in \h{\U}_{1026}$, then, as in the proof of Proposition \ref{prop zeta support}, we write $\chi = \chi_2\chi_3\chi_{19}$ where $\chi_p \in \h{\U}_{p^{m_p}}$. Observe that $\U_{1026} \cong \ZZ/(18)^2$ is generated by $191$ and $325$. Then for $(a,b) \in \ZZ/(18)^2$ we have
\[
    \chi(191^a 325^b) = \chi_3(2)^a \chi_{19}(2)^b.
\]
Thus $v_3(c_\chi) \leq 1$ if and only if $\chi_3(2) = \pm 1$. Identifying $\U_{1026}$ with the dual of $\h{\U}_{1026}$ we choose coordinates so that $\rho(191) = x$ and $\rho(325) = y$. Then in additive coordinates it follows that $\Sigma_{3,1026}$ is the complement of the hyperplane $\H_\Sigma : 2x = 0$. 

Suppose $d$ is coprime to $1026$ such that $M_{6d}(\zeta_{1026}) = 0$. Therefore, in the notation of Theorem \ref{thm unlikely neck}, we have $e = 2$ and $f = 3$. Our computation above together with Theorem \ref{thm unlikely neck} imply that $M_{6d}(\zeta_{1026}) = 0$ if and only if
\[
    \h{\U}_{1026} \subseteq \bigcup_{p\mid d}\H_p \cup \H_2^a \cup \H_\Sigma,
\]
where $\H_2^a$ is the affine hyperplane of all $\chi$ such that $\chi(2) = -1$. Since $2 \equiv 191\cdot 325 \bmod 1026$, the affine hyperplane $\H_2^a$ may be expressed additively as $x + y = 9$. Therefore there is a correspondence between $d$ coprime to $1026$ such that $M_{6d}(\zeta_{1026}) = 0$ and arrangements of lines in $\ZZ/(18)^2$ which, together with $2x = 0$ and $x + y = 9$, cover all of $\ZZ/(18)^2$. One example is illustrated below.
\begin{center}
    \includegraphics[scale=.18]{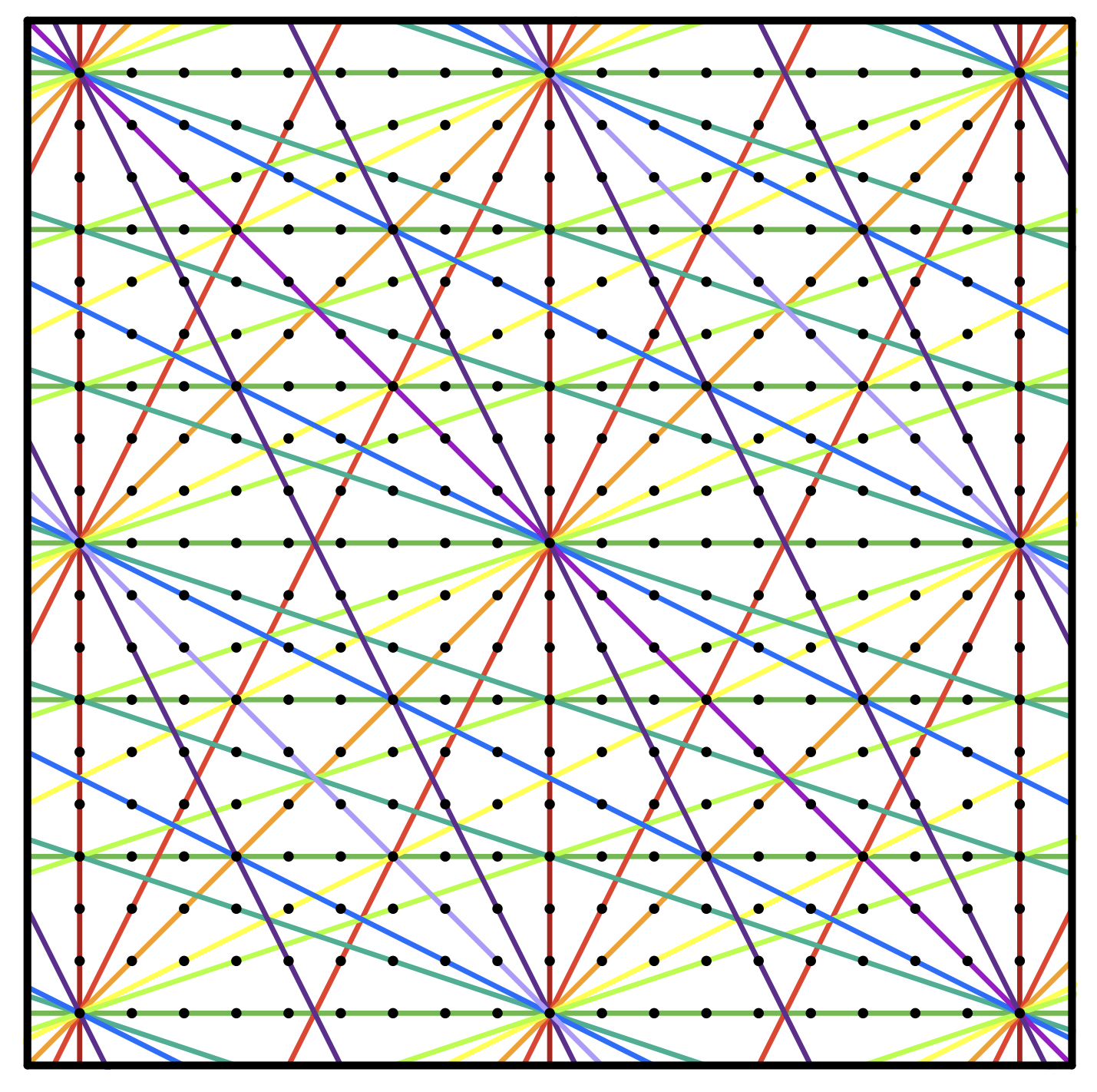}
\end{center}
Such arrangements can be found by starting with $2x = 0$ and $x + y = 9$ and drawing lines until every point is covered. To convert this picture into a concrete solution, we use our choice of coordinates to translate each linear form into a residue modulo 1026 and then pick a prime in that congruence class. For example, the line $2x + 6y = 0$ corresponds to the residue class $191^2\cdot 325^6 \equiv 463 \bmod 1026$. Since $463$ is prime, the line $2x + 6y = 0$ is $\H_{463}$. Applying this procedure to the above diagram we find
\[
    d = 61\cdot 139\cdot  463\cdot  733\cdot  859\cdot  919\cdot  1327\cdot  2797\cdot  2797\cdot  3593 = 84732759227967517764591359639.
\]
Thus the 1026th cyclotomic polynomial divides the 508396555367805106587548157834th necklace polynomial and this is a reflection of the fact that the family of lines depicted above covers $\ZZ/(18)^2$.
\end{example}


\section{Cyclotomic factors of \texorpdfstring{$\Phi_d(x) -1$}{Phid(x) - 1}}
\label{subsec cyclo evals}

In this section we characterize the pairs $(d,m)$ for which $\Phi_d(\zeta_m) = 1$ in terms of hyperplane arrangements covering explicit subsets of $\h{\U}_m$. The structure of this section parallels that of Section \ref{sec neck evals}.

We will make use of the following functions of a real variable $x$ with $d\geq 1$ defined by
\[
    \zeta^x := \exp(2\pi i x) \hspace{.4in} \ve(x) := 2|\sin(\pi x)| \hspace{.4in}  \varphi_d\lfloor x\rfloor := \sum_{e\mid d}\mu(d/e)\lfloor ex \rfloor \equiv \sum_{e\mid d}\lfloor ex\rfloor \bmod 2.
\]
Thus $\zeta^{k/m} = \zeta_m^k$ and $\ve(x)$ is periodic with period $1$ and positive for all non-integral $x$. If $d \in \NN$, let $\ve(x)^{[d]} := \ve(dx)$.

\begin{lemma}
\label{lem polar decomp}
Let $d > 1$ and let $x$ be a real variable, then
\begin{enumerate}
    \item $\displaystyle{\zeta^x - 1 = i(-1)^{\lfloor x \rfloor}\zeta^{x/2} \ve(x).}$\\
    
    \item $\displaystyle{\Phi_d(\zeta^x) = (-1)^{\varphi_d\lfloor x\rfloor}\zeta^{\varphi(d)x/2}\ve(x)^{\varphi_d}.}$
\end{enumerate}
\end{lemma}

\begin{proof}
(1) Recall that $2\sin(\pi x) = -i(\zeta^{x/2} - \zeta^{-x/2})$. Thus,
\[
    \zeta^x - 1 = \zeta^{x/2}(\zeta^{x/2} - \zeta^{-x/2})
    = i\zeta^{x/2}(2\sin(\pi x))
    = i(-1)^{\lfloor x \rfloor} \zeta^{x/2} \big((-1)^{\lfloor x \rfloor}2\sin(\pi x)\big).
\]
The functions $(-1)^{\lfloor x \rfloor}$ and $2\sin(\pi x)$ are both periodic with period $2$. Since
\begin{align*}
    2\sin(\pi(x + 1)) &= -2\sin(\pi x)\\
    (-1)^{\lfloor x + 1 \rfloor} &= - (-1)^{\lfloor x \rfloor},
\end{align*}
it follows that their product has period 1 and
\[
    (-1)^{\lfloor x \rfloor}2\sin(\pi x) = 2|\sin(\pi x)| = \ve(x).
\]
Therefore,
\[
    \zeta^x - 1 = i(-1)^{\lfloor x \rfloor}\zeta^{x/2} \ve(x).
\]

(2) We compute,
\begin{align*}
    \Phi_d(\zeta^x) &= \prod_{e\mid d}(\zeta^{ex} - 1)^{\mu(d/e)}\\
    &= \prod_{e\mid d} (i(-1)^{\lfloor ex \rfloor}\zeta^{ex/2} \ve(ex))^{\mu(d/e)}\\
    &= (-1)^{\sum_{e\mid d}\mu(d/e)\lfloor ex \rfloor} \zeta^{\sum_{e\mid d}\mu(d/e)ex/2} \prod_{e\mid d}\ve(ex)^{\mu(d/e)}\\
    &=(-1)^{\varphi_d\lfloor x\rfloor} \zeta^{\varphi(d)x/2} \ve(x)^{\varphi_d}.
\end{align*}
Note that the factor of $i$ cancels in the third equality because $\sum_{e\mid d} \mu(d/e) = 0$ for $d > 1$.
\end{proof}

Evaluating Lemma \ref{lem polar decomp}.2 at $x = 1/m$ gives us  $\Phi_d(\zeta_m) = (-1)^{\varphi_d\lfloor 1/m\rfloor} \zeta_{2m}^{\varphi(d)} \ve(1/m)^{\varphi_d}$. Therefore $\Phi_d(\zeta_m) = 1$ is equivalent to the following two identities holding simultaneously,
\begin{align}
    \label{eqn phase}(-1)^{\varphi_d\lfloor 1/m\rfloor} \zeta_{2m}^{\varphi(d)} &= 1\\ \label{eqn radial}\ve(1/m)^{\varphi_d} &= 1.
\end{align}
The first equation \eqref{eqn phase} is a purely arithmetic condition. The second equation \eqref{eqn radial} requires us to analyze the support of $\ve(1/m)$ in the $\ZZ[\s\U_m]$-module $\QQ(\zeta_m)^\times$.

For the convenience of additive notation, let 
\[
    \ell(x) := \log \ve(x).
\]
We define the \emph{$m$th cyclotomic module} $\C_m$ to be the $\QQ[\s\U_m]$-module spanned by $\ell(a/m)$ for $a \not\equiv 0 \bmod m$. Let $[q]\ell(a/m) := \ell(qa/m)$ for $q \in \U_m$, then
\[
    \varphi_d\ell(1/m) = \log |\Phi_d(\zeta_m)| = \log \ve(1/m)^{\varphi_d}.
\]
Since
\[
    \ell(-x) = \log |\zeta^{-x} - 1| = \log |\zeta^{-x}(1 - \zeta^x)| = \ell(x),
\]
the action of $\U_m$ on $\C_m$ factors through $\U_m/\langle -1 \rangle$.

Bass \cite[Thm. 2]{Bass} determined the structure of $\C_m$ as a $\QQ[\s\U_m]$-module. The proof of Theorem \ref{thm Bass} is a combination of a Galois equivariant version of the Dirichlet unit theorem and the fact that the cyclotomic units have finite index in the units $\ZZ[\zeta_m]^\times$.

\begin{thm}[Bass]
\label{thm Bass}
Let $m\geq 1$, let $\omega(m)$ denote the number of distinct prime factors of $m$, and let $\mathbf{1}$ denote the trivial representation of $\U_m$. Then
\[
    \C_m \cong \QQ[\s\U_m/\langle -1\rangle] \oplus \mathbf{1}^{\omega(m)-1}.
\]
Therefore the support of $\C_m$ is $\{\chi \in \h{\U}_m : \chi(-1) = 1\} = \H_{-1}$.
\end{thm}

\begin{remark}
Our definition of the $m$th cyclotomic module varies slightly from how Bass defines it. Bass' ($\QQ$-linearized) cyclotomic module $\C_m'$ is defined as the $\QQ$-extension of scalars of the abelian group multiplicatively spanned by $\zeta_m^a - 1$ with $a \not\equiv 0 \bmod m$. There is a natural surjective map $\C_m' \rightarrow \C_m$ given by $\zeta_m^a - 1 \mapsto \log|\zeta_m^a - 1| = \ell(a/m)$ which we claim is an isomorphism. It suffices to show that if $u = \prod_{a = 1}^{m-1}(\zeta_m^a - 1)^{b_a} \in \C_m'$ with $b_a \in \ZZ$ has absolute value 1, then $u$ is a root of unity. Since $\overline{(\zeta_m^a - 1)} = -\zeta_m^{-a}(\zeta_m^a - 1)$, we have
\[
    1 = |u| = u\overline{u} = \zeta u^2,
\]
for some root of unity $\zeta$. Hence $u$ is a square root of a root of unity, and thus is itself a root of unity. Therefore $\C_m \cong \C_m'$.
\end{remark}

The following lemma establishes several useful relations in $\C_m$.

\begin{lemma}
\label{lem rel norm}
Let $m > 1$.
\begin{enumerate}
    \item If $q > 1$ is a natural number not divisible by $m$, then
    \[
        ([q] - 1)\ell(1/m) := \ell(q/m) - \ell(1/m) = \sum_{b=1}^{q-1} \ell(1/m + b/q).
    \]
    
    \item Let $p$ be a prime and suppose that $q = p^e$ divides $m$.
    \begin{enumerate}
        \item If $d < e$, then $\ell(p^d/m) \in \QQ[\s\U_m]\ell(1/m)$.
        \item If $q$ is the largest power of $p$ dividing $m$ and $n = m/q$, then $\varphi_p\ell(1/n) \in \QQ[\s\U_m]\ell(1/m)$.
    \end{enumerate}
    
\end{enumerate}
\end{lemma}

\begin{proof}
(1) Observe that
\[
    |\zeta^{qx} - 1| = \prod_{b=0}^{q-1}|\zeta^{x + b/q} - 1|.
\]
Evaluating at $x = 1/m$ and taking logarithms (which we can because $m \nmid q$) we find
\[
    \ell(q/m) = \sum_{b=0}^{q - 1}\ell(1/m + b/q).
\]

(2a) Part (1) implies that
\[
    \ell(p^d/m) = \sum_{b=0}^{p^d -1}\ell(1/m + b/p^d) = \sum_{b=0}^{p^d -1}\ell\Big(\frac{1 + b(m/p^d)}{m}\Big).
\]
Since $d < e$, we see that $m/p^d$ is divisible by $p$. Hence $1 + b(m/p^d)$ is a unit modulo $m$. Thus $\ell(p^d/m) \in \QQ[\s\U_m]\ell(1/m)$.

(2b) Let $n := m/q$, so that $n$ is coprime to $p$ by assumption. Applying (1) we have
\[
    \varphi_p\ell(1/n) = \sum_{b=1}^{p-1}\ell(1/n + b/p) = \sum_{b=1}^{p-1}\ell\Big(\frac{p + bn}{np}\Big).
\]
Since $n$ and $p$ are coprime and $b$ is a unit modulo $p$, it follows that $p + bn$ is a unit modulo $n$ and modulo $p$, hence a unit modulo $np$. Therefore
\[
    \varphi_p\ell(1/n) \in \QQ[\s\U_m]\ell(1/np) \subseteq \QQ[\s\U_m]\ell(1/m),
\]
where the last inclusion is a consequence of part (2a).
\end{proof}

\begin{prop}
\label{prop cyclo support}
Let $m > 1$ be an integer and let $\chi \in \h{\U}_m$ be a character of modulus $m$. Then $\chi$ is in the support of $\ell(1/m)$ if and only if
\[
    \chi \in \H_{-1} \setminus \bigcup_{p\mid m} \H_p = \{\chi \in \h{\U}_m : \chi(-1) = 1 \text{ and } \chi(p) \neq 1 \text{ for all primes } p \mid m\}.
\]
\end{prop}

\begin{remark}
Recall that by our convention on extending the domains of Dirichlet characters (see Caution \ref{caution},) if $p$ is a prime dividing $m$, then $\chi(p)$ has a well-defined, nonzero value if the conductor of $\chi$ is not divisible by $p$, and otherwise $\chi(p)=0$. If $q$ is the largest power of $p$ dividing $m$, then $\H_p \subseteq \h{\U}_{m/q} \subseteq \h{\U}_m$.
\end{remark}

\begin{proof}
Let $m = q_1q_2\cdots q_k$ be the factorization of $m$ into prime powers where $q_i$ is a power of the prime $p_i$. If $J \subseteq \{1, 2, \ldots, k\}$ is a subset, let $m_J := \prod_{j \in J} q_j$ and let $n_J := m/m_J$.
Lemma \ref{lem rel norm}.2b implies that for each proper subset $J \subset \{1,2,\ldots, k\}$, 
\[
    \Big(\prod_{p\mid m_J}\varphi_{p}\Big)\ell(1/n_J) \in \QQ[\s\U_m]\ell(1/m).
\]
Let $\wt{\Sigma}_m$ denote the support of $\ell(1/m)$, then the support of the above element is $\wt{\Sigma}_{n_J} \setminus \bigcup_{p\mid m_J}\H_{p}$. Hence
\[
    \wt{\Sigma}_m \supseteq \wt{\Sigma}_{n_J} \setminus \bigcup_{p\mid m_J}\H_{p}.
\]
Lemma \ref{lem rel norm}.2 shows that $\C_m$ is generated as a $\QQ[\s\U_m]$-module by $\ell(1/n_J)$ as $J$ ranges over all proper subsets of $\{1,2,\ldots, k\}$ and Theorem \ref{thm Bass} shows that the support of $\C_m$ is $\H_{-1}$. Thus
\[
    \H_{-1} = \bigcup_{J}\wt{\Sigma}_{n_J}.
\]
Therefore
\[
    \wt{\Sigma}_m \supseteq \bigcup_{J} ( \wt{\Sigma}_{n_J} \setminus \bigcup_{p \mid m_J}\H_{p}) \supseteq \H_{-1} \setminus \bigcup_{p\mid m}\H_{p}.
\]
Now we show the reverse inclusion. Lemma \ref{lem rel norm} implies that for each $i$,
\[
    \varphi_{p_i}\ell(1/n_i)
    = \Big(\sum_{b=1}^{p_i-1}[p_i + bn_i]\Big) \ell(1/n_ip_i)
    = \Big(\sum_{b=1}^{p_i-1}[p_i + bn_i]\Big)[q_i/p_i] \ell(1/m).
\]
If $\chi \in \H_{p_i}$, then by definition $\chi$ must have modulus $n_i$ and $\chi(p_i) = 1$. Thus applying the idempotent $e_\chi$ to the right hand side of the above identity we find
\[
    e_\chi \Big(\sum_{b=1}^{p_i-1}[p_i + bn_i]\Big)[q_i/p_i] \ell(1/m) = \sum_{b=1}^{p_i - 1}\chi(p_i + bn_i)\chi(q_i/p_i)\ell(1/m)_\chi = (p_i - 1)\ell(1/m)_\chi.
\]
On the other hand,
\[
    e_\chi \varphi_{p_i}\ell(1/n_i) = (\chi(p_i) - 1) \ell(1/n_i)_\chi = 0.
\]
Therefore $\ell(1/m)_\chi = 0$, which is equivalent to saying that $\chi$ does not belong to $\wt{\Sigma}_m$. Hence
\[
    \wt{\Sigma}_m \subseteq \H_{-1} \setminus \bigcup_{p\mid m}\H_p.\qedhere
\]
\end{proof}

We now prove the main result of this section.

\begin{thm}
\label{thm unlikely cyclo}
Let $d, e, f \geq 1$ and $m > 1$ be integers, let $m'$ be the squarefree part of $m$, and let $\H_3^a\subseteq \h{\U}_m$ be the affine hyperplane $\H_3^a := \{\chi \in \h{\U}_m : \chi(3) = -1\}$. Suppose that
\begin{multicols}{2}
\begin{enumerate}[label = (\roman*)]
    \item $m$ does not divide $def$,
    \item $def$ is squarefree,
    \item $d$ is coprime to $m$,
    \item $e$ divides $m'$,
    \item $f$ divides $m/m'$.
\end{enumerate} 
\end{multicols}
Then $\Phi_{def}(\zeta_m) = 1$ if and only if
\begin{enumerate}
    \item $\displaystyle{\H_{-1} \subseteq \begin{cases} \bigcup_{p\mid md/e} \H_p & \text{if }3 \nmid e\\ \bigcup_{p\mid md/e} \H_p\cup \H_3^a & \text{if }3 \mid e, \end{cases}}$\\
    \item $m$ divides $\varphi(def)$, and\\
    \item $\displaystyle{\sum_{a\mid def}\lfloor a/m \rfloor \equiv \frac{\varphi(def)}{m} \bmod 2}$.
\end{enumerate}
\end{thm}

\begin{remark}
Recall that by Remark \ref{rem sqfree} we lose no generality in assuming that $def$ is squarefree.
\end{remark}

\begin{proof}
As we observed following Lemma \ref{lem polar decomp}, $\Phi_{def}(\zeta_m) = 1$ is equivalent to the triviality of both the phase \eqref{eqn phase} and the radial \eqref{eqn radial} components of $\Phi_{def}(\zeta_m)$. Suppose that the phase component of $\Phi_{def}(\zeta_m)$ is trivial,
\[
    (-1)^{\lfloor\varphi\rfloor_{def}(1/m)} \zeta_{2m}^{\varphi(def)} = 1.
\]
Thus $\zeta_{2m}^{\varphi(def)} = \pm 1$, which is equivalent to $m$ dividing $\varphi(def)$. If $m$ does divide $\varphi(def)$, then comparing exponents of $-1$ in the above identity we conclude that
\[
    \sum_{a\mid def}\lfloor a/m \rfloor \equiv \frac{\varphi(def)}{m} \bmod 2.
\]

Triviality of the radial component of $\Phi_{def}(\zeta_m)$ is equivalent to
\[
    \varphi_{def}\ell(1/m) = 0.
\]
Following the same strategy as Theorem \ref{thm unlikely neck}, we determine the support of $\varphi_{ef}\ell(1/m)$ and then appeal to Theorem \ref{thm key}.

Let $\chi \in \H_{-1} \subseteq \h{\U}_m$ be a character. If $\wt{e}$ is the product of all primes $p$ dividing $e$ such that $\chi(p) = 1$ and $\wt{n} := m/\wt{e}$, then we claim that there is some nonzero constant $c$ such that
\begin{equation}
\label{eqn ef claim}
    (\varphi_{ef}\,\ell(1/m))_\chi = \begin{cases} 0 & \text{if } 3 \mid e \text{ and } \chi \in \H_3^a,\\ c\,\ell(1/\wt{n})_\chi & \text{otherwise.}  \end{cases}
\end{equation}

First we finish the proof supposing that we have shown \eqref{eqn ef claim}. Proposition \ref{prop cyclo support} implies that $\ell(1/\wt{n})_\chi = 0$ if and only if $\chi(p) = 1$ for some prime $p \mid \wt{n}$, and any such prime must divide the factor $m/e$ of $\wt{n}$ by the definition of $\wt{n}$. Therefore, the support $\Sigma$ of $\varphi_{ef}\,\ell(1/m)$ is
\[
    \Sigma = \begin{cases}\H_{-1} \setminus \bigcup_{p\mid m/e}\H_p & \text{if }3 \nmid e,\\ \H_{-1} \setminus (\bigcup_{p\mid m/e}\H_p \cup \H_3^a) & \text{if } 3 \mid e. \end{cases}
\]
Thus Theorem \ref{thm key} implies $\varphi_{def}\ell(1/m) = 0$ if and only if 
\[
    \H_{-1} \subseteq \begin{cases} \bigcup_{p\mid dm/e} \H_p & \text{if }3 \nmid e\\ \bigcup_{p\mid dm/e} \H_p\cup \H_3^a & \text{if }3 \mid e. \end{cases}
\]

All that remains is to prove \eqref{eqn ef claim}. We use the factorization $\varphi_{ef} = \prod_{p\mid ef} \varphi_p$ (which uses \textit{(ii)}) to analyze $(\varphi_{ef}\,\ell(1/m))_\chi$ one prime at a time.

Let $p$ be a prime dividing $f$ and let $n := m/p$. Then $p$ divides $n$ by \textit{(v)}. Lemma \ref{lem rel norm}.1 implies that
\[
    \varphi_p\, \ell(1/m) = \sum_{k\in\s\U_p}[1 + kn]\ell(1/m).
\]
If $\chi$ has conductor dividing $n$, then
\[
    (\varphi_p\,\ell(1/m))_\chi = \sum_{k\in\s\U_p}\chi(1 + kn)\ell(1/m)_\chi = (p - 1)\ell(1/m)_\chi.
\]
If the conductor of $\chi$ does not divide $n$, then write $\chi = \chi_n \chi_p$ where $\chi_n$ has conductor dividing $n$ and $\chi_p$ has conductor $p^{m_p}$ with $m_p = v_p(m) > 1$ (the inequality uses \textit{(v)}.) Thus
\[
    (\varphi_p\,\ell(1/m))_\chi = \sum_{k\in\s\U_p}\chi(1 + kn)\ell(1/m)_\chi
    = \sum_{k\in\s\U_p} \chi_p(1 + kp^{m_p-1}) \ell(1/m)_\chi
    = -\ell(1/m)_\chi,
\]
where the last equality follows from the observation that $\chi_p(1 + k p^{m_p-1})$ ranges over the the non-trivial $p$th roots of unity as $k$ ranges over $\U_p$. Hence in either case there is some nonzero constant $c$ such that
\begin{equation}
\label{eqn f}
    (\varphi_p\,\ell(1/m))_\chi = c\s\ell(1/m)_\chi. 
\end{equation}

Next let $p$ be a prime dividing $e$,
so that $n := m/p$ is coprime to $p$ by \textit{(iv)} and $n > 1$ by \textit{(i)}. Observe that
\[
    (\varphi_p\,\ell(1/m))_\chi = \ell(1/n)_\chi - \ell(1/m)_\chi.
\]
If the conductor of $\chi$ does not divide $n$, then
\[
    (\varphi_p\,\ell(1/m))_\chi = -\ell(1/m)_\chi.
\]
Suppose the conductor of $\chi$ divides $n$. Lemma \ref{lem rel norm}.1 implies that
\begin{equation}
\label{eqn 1}
    (\chi(p) - 1)\ell(1/n)_\chi = (\varphi_p\,\ell(1/n))_\chi 
    = \sum_{k \in \s\U_p} \chi(p + kn) \ell(1/m)_\chi 
    = \chi(p)(p - 1)\ell(1/m)_\chi.
\end{equation}

If $\chi(p) \neq 1$, then
\[
    (\varphi_p\,\ell(1/m))_\chi = \ell(1/n)_\chi - \ell(1/m)_\chi = \Big(\frac{p-1}{1 - \overline{\chi(p)}} - 1\Big)\ell(1/m)_\chi.
\]
The coefficient of $\ell(1/m)_\chi$ vanishes if and only if $\chi(p) = 1/(2 - p)$. Since $\chi(p)$ is a root of unity, it must be the case that $p = 3$ and $\chi(3) = -1$.

If $\chi(p) = 1$, then $\ell(1/m)_\chi = 0$ by \eqref{eqn 1} and thus 
\[
    (\varphi_p\,\ell(1/m))_\chi = \ell(1/n)_\chi.
\]

Hence if $p \mid e$ and $n = m/p$, then
\begin{equation}
\label{eqn e}
    (\varphi_p \ell(1/m))_\chi = \begin{cases} \ell(1/n)_\chi & \text{if $\chi \in \H_p$,}\\ 0 & \text{if $p = 3$ and $\chi \in \H_3^a$}\\ c\s\ell(1/m)_\chi & \text{otherwise, for some nonzero constant $c$.} \end{cases}
\end{equation}
Together \eqref{eqn f} and \eqref{eqn e} prove our claim \eqref{eqn ef claim}.
\end{proof}

\begin{example}
Let $m = 24$ and suppose we want to find an integer $d$ coprime to $24$ such that $\Phi_{3d}(\zeta_{24}) = 1$. The group of Dirichlet characters $\h{\U}_{24}$ is a 3 dimensional $\FF_2$-vector space. Let $\rho: \U_{24} \rightarrow \h{\FF}_2^3$ be the choice of coordinates such that
\[
    \rho(13) = x, \hspace{.35in} \rho(17) = y, \hspace{.35in} \rho(19) = z.
\]
Following Theorem \ref{thm unlikely cyclo}, we begin by looking for some $d$ such that the plane $\H_{-1}: x + y + z = 0$ is covered by $\H_2$, $\H_3^a$, and the union of the $\H_p$ with $p \mid d$. Since $13 \equiv 19 \equiv 1 \bmod 3$ and $17 \equiv 1 \bmod 8$, it follows that $\h{\U}_3$ is the subspace $x = z = 0$ and $\h{\U}_8$ is the subspace $y = 0$. Then $\H_2 \subseteq \h{\U}_3$ consists of the single point $(0,0,0)$ and $\H_3^a\subseteq \h{\U}_8$ is the subspace $x = 1$, $y = 0$, which intersects $x + y + z = 0$ at the point $(1,0,1)$. Therefore it suffices for $\bigcup_{p\mid d}\H_p$ to cover the two points $(1,1,0)$ and $(0,1,1)$. For example, the lines $\H_{5} = \H_{13\cdot 17} : x + y = 0$ and $\H_{11} = \H_{17\cdot 19} : y + z = 0$ suffice.

If $d = 55 =  5\cdot 11$, then $\varphi(3d) = 80$ is not divisible by $24$, hence $\Phi_{3\cdot 5\cdot 11}(\zeta_{24}) \neq 1$ by Theorem \ref{thm unlikely cyclo}. On the other hand, if $d = 385 = 5 \cdot 7 \cdot 11$, then  $\varphi(3\cdot 5\cdot 7\cdot 11) = 480$ is divisible by 24 and
\[
    \sum_{a \mid 3d} \lfloor a/m \rfloor = 90 \equiv 20 =\frac{\varphi(3d)}{m} \bmod 2.
\]
Thus Theorem \ref{thm unlikely cyclo} implies that $\Phi_{3d}(\zeta_m) = \Phi_{1155}(\zeta_{24}) = 1$.
\end{example}

\end{document}